\documentclass[leqno]{report}
\usepackage{amsmath,amsthm,amsfonts,amssymb,setspace,graphicx,
diagrams,float, floatflt,wrapfig,verbatim}
\newcommand{\frontmatter}{\clearpage \pagenumbering{roman}}
\newcommand{\mainmatter}{\clearpage \pagenumbering{arabic}}

\usepackage{graphicx}
\usepackage{amssymb}
\usepackage{relsize}
\usepackage{amsmath,amscd}
\usepackage[mathscr]{eucal}
\newtheorem{thm}{\hspace{1cm}Theorem}[section]

\newtheorem{lem}[thm]{\hspace{1cm}Lemma}
\newtheorem{sublem}[thm]{\hspace{1cm}Sublemma}
\newtheorem{prop}[thm]{\hspace{1cm}Proposition}

\theoremstyle{definition}
\newtheorem{defn}[thm]{\hspace{1cm}Definition}
\theoremstyle{remark}
\newtheorem{rem}[thm]{\hspace{1cm}Remark}
\numberwithin{equation}{section}

\newenvironment{myproof}[1]{\emph{#1}.  }{\qed\vspace{.3cm}}

 \newcommand{\s}{\mathcal{S}}

 \newcommand{\Real}{\mathbb{R}}
 \newcommand{\Natural}{\mathbb{N}}
 
 \newcommand{\Complex}{\mathbb{C}}
 \newcommand{\Rational}{\mathbb{Q}}
 
 \newcommand{\Int}{\mathbb{Z}}

 \newcommand{\Gh}{\mathfrak{h}}

 \newcommand{\GL}[2]{\mathrm{GL}_{#1}(#2)}
 
 \newcommand{\SL}[2]{\mathrm{SL}_{#1}(#2)}

 \newcommand{\SO}[2]{\mathrm{SO}_{#1}(#2)}
  
 \newcommand{\SOR}[1]{\mathrm{SO}(#1)}
 \newcommand{\SOM}[2]{\mathrm{SO}(#1,#2)}

 \newcommand{\intd}{\,\mathrm{d}}

 \newcommand{\half}{\frac{1}{2}}

 \newcommand{\vectj}{\mathbf{j}}
 \newcommand{\vecti}{\mathbf{i}}

 \newcommand{\conj}{\mathbf{c}}

 \newcommand{\Id}{\mathrm{Id}}

 \newcommand{\Go}{\mathfrak{o}}

 \newcommand{\inv}{^{-1}}
 \newcommand{\scrF}{\mathscr{F}}
  \newcommand{\scrG}{\mathscr{G}}

 \newcommand{\scrS}{\mathscr{S}}

 \newcommand{\scrC}{\mathscr{C}}

 \newcommand{\Spos}{\mathrm{\Spos}}

 \newcommand{\Rept}{\mathrm{Re}}
 \newcommand{\Impt}{\mathrm{Im}}
 \newcommand{\bbH}{\mathbb{H}}
 \newcommand{\bbS}{\mathbb{S}}

 \newcommand{\coh}{\mathscr{C}_{\mathbf{H}}}
 \newcommand{\intrr}{\mathrm{Int}}
 \newcommand{\finishdisp}{\vspace{0.3cm} \]}
 \newcommand{\startdisp}{\vspace{0.3cm}\[ }
 \newcommand{\starteqn}{\vspace{0.3cm}\begin{equation}}
 \newcommand{\finisheqn}{\vspace{0.3cm}\end{equation}}
 \newcommand{\starteqnarraynn}{\vspace*{0.1cm}\begin{eqnarray*}}
 \newcommand{\finisheqnarraynn}{\vspace*{2cm}\end{eqnarray*}}
  \newcommand{\starteqnarray}{\vspace*{0.1cm}\begin{eqnarray}}
 \newcommand{\finisheqnarray}{\vspace*{2cm}\end{eqnarray}}
 
 \newcommand{\Fix}{\mathrm{Fix}}

 \newcommand{\BDelta}{\mathbf{\Delta}}
 
 \newcommand{\forallindisp}{\quad\text{for all}\;}

 \newcommand{\bmth}{\boldmath}
 \newcommand{\ubmth}{\unboldmath}

  \newcommand{\red}{\mathrm{red}}
  \newcommand{\MtwoN}{\mathrm{M}_2^{\mathrm{N}}}

  \newcommand{\scrK}{\mathscr{K}}
  \newcommand{\Iso}{\mathrm{Iso}}

\begin{document}
\pagestyle{headings}
\frontmatter
\begin{abstract}
\begin{center}
A fundamental domain
of Ford type for
$\SO{3}{\Int[\vecti]}\backslash \SO{3}{\Complex}/\SOR{3}$,
and for $\SOM{2}{1}_{\Int}\backslash\SOM{2}{1}/\SOR{2}$.\\
May  2006
\end{center}
Let $G=\SO{3}{\Complex}$,
$\Gamma=\SO{3}{\Int[\vecti]}$, $K=\SOR{3}$, and let
$X$ be the locally symmetric space $\Gamma\backslash G/K$.
In this paper, we write down explicit equations defining a fundamental
domain for the action of $\Gamma$ on $G/K$.  The fundamental
domain is well-adapted for studying the theory of $\Gamma$-invariant
functions on $G/K$.  We write down equations defining a fundamental domain
for the subgroup $\Gamma_{\Int}=\SOM{2}{1}_{\Int}$ of
$\Gamma$ acting on the symmetric space $G_{\Real}/K_{\Real}$,
where $G_{\Real}$ is the split real form $\SOM{2}{1}$ of $G$
and $K_{\Real}$ is its maximal compact subgroup $\SOR{2}$.
We formulate a simple geometric relation between the fundamental
domains of $\Gamma$ and $\Gamma_{\Int}$ so described.  These
fundamental domains are geared towards the detailed study of the
spectral theory of $X$ and the embedded subspace
 $X_{\Real}=\Gamma_{\Int}\backslash G_{\Real}/K_{\Real}$.
\end{abstract}
\begin{titlepage}
\begin{singlespace}
\vspace*{50pt}
\begin{center} \Huge \bfseries
A fundamental domain
of Ford type for
$\SO{3}{\Int[\vecti]}\backslash \SO{3}{\Complex}/\SOR{3}$,
and for $\SOM{2}{1}_{\Int}\backslash\SOM{2}{1}/\SOR{2}$.
\end{center}
\vspace{2.25in}
\begin{center}
A Paper
\end{center}
\vspace{2.25in}
\begin{center}
by \\ Eliot Brenner
\end{center}
\begin{center}
Center for Advanced Studies in Mathematics\\
at Ben Gurion University
\end{center}
\begin{center}
May 2006
\end{center}
\end{singlespace}
\end{titlepage}
\clearpage
\thispagestyle{empty}
\begin{center}
\large\vspace*{1in}
ACKNOWLEDGMENTS
\end{center}
\vspace{0.25in}
The Institute for the Advanced Study of Mathematics
at Ben-Gurion University provided support
and a pleasant and stimulating working environment
during the writing of this paper.  The author thanks
Mr. Tony Petrello for additional financial assistance.
\mainmatter
\section{Introduction}
The author has undertaken, in Chapter 1
of \cite{brennerthesis}, a generalization of the classical
theory of Ford fundamental domains (see \S2.2 of \cite{iwaniec95})
for Fuchsian groups
to a wide class of group actions
including, in particular,
$\Gamma_n=\SL{n}{\Int[\vecti]}$ acting on $G_n=\SL{n}{\Complex}/
\mathrm{SU}(n)$
and $\mathrm{GL}(n,\Int)$ acting on $\mathrm{GL}(n,\Real)/\mathrm{SO}(n)$.
In the latter case, the fundamental domains obtained
coincide  with the $F_n$ studied by
D. Grenier in \cite{grenier88} and
\cite{grenier93} (allowing for the isomorphism of the symmetric space
$G/K$ with the quadratic model $P$).  For this reason,
we adopt the terminology \textit{Grenier domains} for the generalized
Ford domains.
A major theme of Grenier's work in these articles is that
the $F_n$ for different $n$ are best considered as part of an inductive
scheme, since $F_m$ for $m<n$ appear both in the definition
of $F_n$ and in his construction of the Satake compactifications
of the locally symmetric space $\mathrm{GL}(n,\Int)\backslash
\mathrm{GL}(n,\Real)/\mathrm{SO}(n)$.  The base case
of Grenier's inductive scheme is (ignoring the center
of $\mathrm{GL}(n,\Real)$) provided by Dirichlet's classical
fundamental domain for $\SL{2}{\Int}$ acting on the upper half
plane.  The results of this paper may be viewed as providing
the base case for an inductive scheme of the same type corresponding
to the sequence of locally symmetric spaces in
\eqref{eqn:ourlocallysymmetricspace},
below.  Note that the base case for this ``orthogonal"
sequence is considerably more complicated than the base
case for Grenier's ``general linear" sequence.

We take advantage of the well-known isomorphism
\startdisp
\SL{2}{\Complex}/\{\pm I\}
\stackrel{\cong}{\longrightarrow} \SO{3}{\Complex},
\finishdisp
specified at the beginning of \S\ref{sec:lattice}, to identify
the lattice $\SO{3}{\Int[\vecti]}$ with a group of fractional
linear transformations acting on $\bbH^3$.  The purpose of the
present paper is to state explicitly what this arithmetic subgroup
is in explicit matrix terms (Proposition
\ref{prop:inversematrixexplicitdescription})
and give an appropriate fundamental domain for the natural
action on hyperbolic 3-space
(Proposition \ref{prop:gammafunddom}).

Proposition \ref{prop:inversematrixexplicitdescription}, below, has
immediate application in the author's ongoing study
(joint with F. Spinu) of a particular generalization
of Selberg's zeta function.  The
\boldmath\textbf{three-dimensional, vector
Selberg zeta function associated
to a Kleinian group $\Gamma$ and
a unitary representation $\chi$ of $\Gamma$}\unboldmath\hspace*{.25mm} was
recently defined by
J.S. Friedman (following Selberg, A.B. Venkov, and others)
by
\starteqn\label{eqn:eulerproduct}
Z_{\Gamma,\chi}(s)=
\prod_{\{\gamma\}}\prod_{k=0}^{\infty}\det(1-\chi(\gamma)
N_0(\gamma)^{-s-k}),\;
\text{for}\; \Rept s\gg 0.
\finisheqn
In the ``Euler product" expression of \eqref{eqn:eulerproduct},
$\{\gamma\}$ ranges over $\Gamma$-conjugacy
classes of primitive hyperbolic elements in $\Gamma$
and $N_0(\gamma)$ denotes the length
of the closed geodesic on
$\Gamma\backslash G/K$ corresponding to $\gamma$.  The meromorphic
continuation of $Z_{\Gamma,\chi}$ (or, more precisely,
of its logarithmic derivative $Z'/Z$) to the entire
complex domain is closely related to an explicit
form of the Selberg trace formula, worked out,
for example, in \cite{friedman05} in parallel
to \cite{egm98}.  It is of obvious interest
to obtain relations between the $Z_{\Gamma,\chi}$
of the members of a pair of lattices $(\Gamma,\Gamma')$,
where $\Gamma$ and $\Gamma'$ are related in various ways.  
For example, in the case of $(\Gamma,\Gamma')$ a pair
of Fuchsian groups, with $\Gamma'\subseteq \Gamma$
and $[\Gamma:\Gamma']<\infty$ (with $Z_{\Gamma,\chi}$
defined similarly for Fuchsian groups), \cite{venkzog82}
gave a formula which is loosely called a ``factorization
formula", because in the case $\Gamma'$ normal in $\Gamma$,
it specializes to a \textit{bona fide}
factorization of $Z_{\Gamma',\chi}$ as the product of
of the $Z_{\Gamma,\chi_i}$, where $\chi_i$ ranges over the
irreducible direct summands
of $\mathrm{Ind}_{\Gamma'}^{\Gamma}\chi$.
In \cite{brennerspinu}, we will consider
such relations for pairs $(\Gamma,\Gamma')$ of commensurable
Kleinian groups in general and in particular, for
the pair
$\left(\conj\inv(\mathrm{SO}_3(\Int[\vecti])),
\mathrm{PSL}_2(\Int[\vecti])\right)$.  It is
clear from the definition \eqref{eqn:eulerproduct}
that one needs to develop concrete understanding
of the relations between the hyperbolic
conjugacy classes of the groups in question,
and Proposition \ref{prop:inversematrixexplicitdescription},
below, lays
the foundations for that study.

In \S\ref{sec.concludingremarksPacific.tex},
we discuss the application
of fundamental domains to the study of
a more general class of spectral zeta functions.

Based on the $\mathrm{SL}_n/\mathrm{GL}_n$ examples in
the literature, one can speculate on future applications
of exact fundamental domains to traditional problems in number theory.  Some diverse examples of applications of Grenier's domain for $\mathrm{GL}_n(\Int)$, acting
on the space of positive-definite real matrices $P_n$, include
the proof in \cite{chjt98} of a bound on the first nontrivial
eigenvalue of the Laplacian for the case $n=3$, the application
in \cite{vulakh04} to the problem of finding a fundamental
system of units in a number field, and most recently
the investigations of \cite{sarnakstromb} into the minima
of Epstein's zeta function.  It seems likely that, as the
detailed study of automorphic functions on quotients of $\SO{n}{\Complex}$
and its real forms becomes more developed, the exact fundamental
domains, which the present paper specifies in the ``base case" $n=2$,
will play a large role in investigating certain zeta functions associated
to these arithmetic quotients.

We mention the relation of Propositions
\ref{prop:inversematrixexplicitdescription}, \ref{prop:gammafunddom},
and \ref{prop:funddomsrelation}, below,
to some results already in the literature.  First, M. Babillot,
at Lemma 3.2 of \cite{babillot02}, constructs a fundamental
domain for $\mathrm{SO}(2,1)_{\Int}$ acting naturally
on the hyperboloid of one sheet.  The method there bypasses
results like Propositions
\ref{prop:inversematrixexplicitdescription} and \ref{prop:gammafunddom} by embedding $\mathrm{SO}(2,1)_{\Int}$ 
as a subgroup of a triangle group of index two.
The fundamental domain so obtained is used to give
a constructive proof that $\mathrm{SO}(2,1)_{\Int}$ acts
with finite covolume, in order that a general theorem
can be applied to solve a lattice-point counting problem.
Also, there is a well-developed theory of \textit{splines},
which are  models for the arithmetic
quotients of $\Rational$-rank-one groups, in a way
different from, but related to, (Grenier) fundamental domains.  
For a recent treatment with a general existence
theorem and references, see \cite{yasakispines}.  It would
be interesting (and possibly useful for cohomology
calculations of the sort undertaken in \cite{yasakisu21})
to determine precisely the relation of ``duality"
that seems to exist between the splines and Grenier
fundamental domains.  However, this is more relevant to higher rank, and therefore, belongs more
to the continuation of the study undertaken in \cite{brennerthesis}
than to the study at hand.
Finally, Chapters 7--9 of \cite{egm98} contain a
treasure-trove of arithmetic-geometric information on the Kleinian
groups $\mathrm{SL}_2(\mathfrak{o}_K)$, where $\mathfrak{o}_K$
denotes the ring of integers in the imaginary quadratic number field $K$.  This paper's treatment of 
$\conj\inv(\mathrm{SO}_3(\Int[\vecti]))$ runs in parallel
to these chapters of \cite{egm98} and provides a foundation
for the future study of automorphic forms on the complex
orthogonal groups in the explicit style of the
subsequent chapters of \cite{egm98}.

The verifications of all the principal propositions
of the present paper are elementary, though lengthy, and they are
not needed for the envisioned applications of the results.
Accordingly,
many details of proofs are omitted and the interested reader
is referred to the electronically archived preprint \cite{brennso3preprint}
for them.

\section{Representation of $\SO{3}{\Int[\vecti]}$ as a lattice
in $\SL{2}{\Complex}$}\label{sec:lattice}
We begin by establishing some basic notational conventions.

Let $n$ be a positive integer and $\Go$ a ring.  We will
use $\mathrm{Mat}_n(\Go)$ to denote the set of all\linebreak
$n$-by-$n$ square matrices with coefficients in $\Go$.
We reserve use the Greek letters $\alpha$, and so on,
for the elements of $\mathrm{Mat}_n(\Go)$,
and the roman letters $a,b,c,d$ and so on,
for the entries of the matrices.
We will denote scalar mutliplication
on $\mathrm{Mat}_n(\Go)$ by simple juxtaposition.
Thus, if $\Go=\Int[\vecti]$,
$\ell\in\Int[\vecti]$ and $\alpha\in\mathrm{Mat}_2(\Int[\vecti])$,
then
\startdisp
\alpha=\begin{pmatrix}a&b\\c&d\end{pmatrix}\;\text{implies}\;
\ell\alpha=\begin{pmatrix}\ell a&\ell b\\\ell c&\ell d\end{pmatrix}.
\finishdisp
The letters $p,q,r,s$ will be reserved
to denote a quadruple of elements of $\Go$
such that $ps-rq=1$.
In what follows, we normally have $\Go=\Int[\vecti]$,
whenever $\alpha$ is written
with entries $p$ through $s$.  Therefore,
\startdisp\alpha=
\begin{pmatrix}p&q\\
r&s\end{pmatrix}\in \SL{2}{\Int[\vecti]},
\finishdisp
unless stated otherwise.

We will denote a conjugation action of a group on a space $V$
by $\conj_V$, when the context makes clear
what this action is.
For example, if $H$ is a linear Lie group and $\Gh$ the Lie
algebra of $H$, then we have
\startdisp
\conj_{\Gh}(h)X=hXh\inv,\forallindisp h\in H,\;X\in\Gh.
\finishdisp
Note that the morphism $\conj_{\Gh}(h)$ is the image under the Lie
functor of the
usual conjugation $\conj_H(h)$ on the group level.
Using $\mathrm{SL}(V)$ to denote the group of unimodular transformations
of a vector space $V$, it is easy to see that
\starteqn\label{eqn:conjalwaysmorphism}
\conj_{\Gh}: H\rightarrow SL(\Gh)\;\text{is a Lie group morphism}.
\finisheqn

Henceforth, whenever $H$ is a group acting on a Lie algebra
$\Gh$ by conjugation, we will omit the subscript $\Gh$.
Thus, we define
\startdisp
\conj:=\conj_{\Gh},
\finishdisp
when we are in the situation of \eqref{eqn:conjalwaysmorphism}.

Except in \S\ref{sec.explicitfd}, we will use the
notation $G=\SO{3}{\Complex}$, $\Gamma=\SO{3}
{\Int[\vecti]}$.
We use $B$ to denote the half-trace form on $\mathfrak{sl}_2(\Complex)$, the Lie
algebra of traceless $2$-by-$2$ matrices.  That is,
\startdisp
B(X,Y)=\half\mathrm{Tr}(XY).
\finishdisp
We use the notation $\beta'=\{X'_1,X'_2,Y'\}$ for the ``standard"
basis of $\mathfrak{sl}_2(\Complex)$, where
\starteqn\label{eqn:betaprime}
\begin{gathered}
X_1'=\begin{pmatrix}0&1\\
0&0
\end{pmatrix},\qquad X_2'=\begin{pmatrix}0&0\\
1&0\end{pmatrix},\\
\text{and}\quad
Y'=\begin{pmatrix}1&0\\
0&-1
\end{pmatrix}
\end{gathered}
\finisheqn

The following properties of $B$ are verified either immediately
from the definition or by straightforward calculations.
\begin{itemize}
\item[\textbf{B1}] $B$ is nondegenerate.
\item[\textbf{B2}]  Setting
\starteqn\label{eqn:beta}
\begin{gathered}
X_1=X_1'+X_2',\qquad X_2=\vecti(X_1'-X_2'),\\
\text{and}\quad Y=Y',
\end{gathered}
\finisheqn
we obtain an orthonormal basis $\beta=\{X_1,X_2,Y\}$, with
respect to the bilinear form $B$.
\item[\textbf{B3}] $B$ is invariant under the conjugation action of $\SL{2}{\Complex}$, meaning
that
\startdisp
B(X,Y)=B(\conj(g)Z,\conj(g)W),\forallindisp Z,W\in\mathfrak{sl}_2(\Complex),\,
g\in\SL{2}{\Complex}.
\finishdisp
\end{itemize}
By \textbf{B3}, $\conj$ is a morphism of $\SL{2}{\Complex}$ into $G$.
The content of part (a) of Proposition \ref{prop:conjsurjection} below
is that the morphism
$\conj$ just described is an epimorphism.

As a consequence of \textbf{B1} and \textbf{B2}, we have that
\starteqn\label{eqn:halftracequadratic}
B(x^1_1 X_1+x^1_2 X_2+y^1 Y,x^2_1 X_1+x^2_2 X_2+y^2Y)=
x^1_1x^2_1+x^1_2x^2_2+y^1y^2,\;
x^i_j, y\in\Complex.
\finisheqn
For any bilinear form $B$ on a vector space $V$,
we use $\mathrm{O}(B)$ to denote the group of linear transformations of $V$
preserving $B$, and we use $\mathrm{SO}(B)$ to denote the unimodular
subgroup of $\mathrm{O}(B)$.
If $B$ is as in \eqref{eqn:halftracequadratic}, then the isomorphism,
\starteqn\label{eqn:concreteorthequiv}
\mathrm{SO}(B)\cong G,
\finisheqn
induced by the identification of the vector space $\mathfrak{sl}_2(\Complex)$
with $\Complex\langle X_1,X_2,Y\rangle$, puts a system of coordinates
on $G$.  Part (b) of Proposition \ref{prop:conjsurjection},
below, will describe the epimorphism $\conj:\SL{2}{\Complex}\rightarrow G$
in terms of these coordinates.
\begin{prop} \label{prop:conjsurjection} With $G$, $\conj$
as above, we have\begin{itemize}
\item[(a)] The map $\conj$ induces an isomorphism
\startdisp
\SL{2}{\Complex}/\{\pm I\}\stackrel{\cong}{\longrightarrow} G
\finishdisp
of Lie groups.
\item[(b)]  Relative to the standard coordinates
on $\SL{2}{\Complex}$ and the coordinates on $G$ induced from
the orthonormal basis $\beta$ of $\mathfrak{sl}_2(\Complex)$, as defined
in \eqref{eqn:beta}, the
epimorphism $\conj: \SL{2}{\Complex}\rightarrow G$
has the following coordinate expression.
\starteqn\label{eqn:imagematrixconj}
\conj\left(\begin{pmatrix}a&b\\c&d\end{pmatrix}\right)=
\begin{pmatrix}\frac{a^2-c^2+d^2-b^2}{2}&
\frac{\vecti(a^2-c^2+b^2-d^2)}{2}&cd-ab\vspace*{0,2cm}\\
\frac{\vecti(b^2+d^2-a^2-c^2)}{2}&\frac{a^2+c^2+b^2+d^2}{2}&\vecti(ab+cd)
\vspace*{0.2cm}\\
-ac+bd&\vecti(ac+bd)&ad+bc\end{pmatrix}.
\finisheqn
\end{itemize}
\end{prop}
We establish some further
notational conventions regarding conjugation mappings.  Whenever
a matrix group $H$ has a conjugation action $\conj_V$ on a \textit{finite dimensional
vector space} $V$ over a field $F$,
each basis $\beta$ of $V$ naturally induces a morphism
\starteqn\label{eqn:conjwithresptobasis}
\conj_{V,\beta}:H\rightarrow\GL{N}{F},\;\text{where}\; N=\dim V.
\finisheqn
Let $\beta$, $\beta'$ be two bases of $V$.
Write $\alpha^{\beta\mapsto\beta'}$ for
the change-of-basis matrix from $\beta$ to $\beta'$.  That is,
if $\beta$, $\beta'$ are written as $N$-entry row-vectors, then
\starteqn\label{eqn:changeofbasis}
\beta\alpha^{\beta\mapsto\beta'}=\beta'.
\finisheqn
Then elementary linear algebra tells us that
\starteqn\label{eqn:conjchangeofbasis}
\begin{aligned}
\conj_{V,\beta}&=&&\conj_{\GL{N}{F}}\left(\left(\alpha^{\beta\mapsto\beta'}\right)\inv\right)\conj_{V,\beta'}\\
&=&&\conj_{\GL{N}{F}}\left(\alpha^{\beta'\mapsto\beta}\right)\conj_{V,\beta'}.
\end{aligned}
\finisheqn
Assuming that $c_V$ is injective, and writing $c_V\inv$ for the
left-inverse of $c_V$, we calculate from \eqref{eqn:conjchangeofbasis}
that
\starteqn\label{eqn:changeofbasisconj}
\conj_{V,\beta}\conj_{V,\beta'}\inv\in\mathrm{Aut}(\GL{N}{F})
\;\text{is given by}\;
\conj_{\GL{N}{F}}\left(\alpha^{\beta\mapsto\beta'}\right).
\finisheqn
In keeping with the practice established after \eqref{eqn:conjalwaysmorphism},
we will omit the subscript $\Gh$ when $H$ is a Lie group
acting on its Lie algebra by conjugation.  Thus, for any basis $\beta$
of $\Gh$,
\startdisp
\conj_{\beta}:=\conj_{\Gh,\beta}.
\finishdisp
Generally speaking, whenever we  fix a single basis $\beta$
for $\Gh$ we will blur the distinction between $\conj$ and $\conj_{\beta}$.
For example, in this paper, whenever $H=\SL{2}{\Complex}$ and $V=\mathrm{Lie}
(H)$, we will write $\conj$ to denote both the ``abstract" morphism $\conj$
of $H$ into $\mathrm{Aut}(V)$ and the linear morphism $\conj_{\beta}$
of $H$ into $\GL{3}{\Complex}$, where $\beta$ is the orthonormal
basis for $\mathrm{Lie}(H)$ defined in \eqref{eqn:beta}.  Whenever
the linear morphism into $\GL{3}{\Complex}$ is induced
by a basis $\beta'\neq\beta$, the notation $\conj_{\beta'}$ will
be used.

We now turn our attention to the description of
the inverse image $\conj\inv(\Gamma)$
as a subset of $\SL{2}{\Complex}/\{\pm I\}$ with respect
to the standard coordinates of $\SL{2}{\Complex}$.  According to Proposition
\ref{prop:conjsurjection}, this amounts to describing the quadruples
\starteqn\label{eqn:quadruplecond}
(a,b,c,d)\in \Complex^4,\;\text{with $ad-bc=1$, and the entries
of the right-side of \eqref{eqn:imagematrixconj} integers.}
\finisheqn
Describing the quadruples meeting conditions
\eqref{eqn:quadruplecond} will
be the subject of the remainder of this section,
culminating in Proposition \ref{prop:inversematrixexplicitdescription}.

\vspace*{0.3cm}\noindent\boldmath\textbf{Conventions regarding multiplicative
structure of $\Int[\vecti]$}.  \unboldmath
Before stating the proposition, we establish
certain conventions we will use when dealing with the multiplicative
properties of the Euclidean ring  $\Int[\vecti]$.  First,
it is well-known that $\Int[\vecti]$ is a Euclidean, hence principal, ring.
That $\Int[\vecti]$ is principal
means that all ideals $\mathscr{I}$ of $\Int[\vecti]$ are generated by a single
element $m\in\Int[\vecti]$, so that every $\mathscr{I}$ is of the form $(m)$.
However, there is an unavoidable ambiguity in the choice of generators
caused by the presence in $\Int[\vecti]$
of four units, $\vecti^{j}$, for $j\in\{0,\ldots, 3\}$, in $\Int[\vecti]$.  We
will adopt the following convention to sidestep the ambiguity caused
by the group of units.
\begin{defn}  We refer to the following subset of $\Complex^{\times}$
as the \textbf{standard subset}
\starteqn\label{eqn:standardset}
\{z\in\Complex^{\times}\;|\; \Rept(z)> 0,\, \Impt(z)\geq 0\}.
\finisheqn
That is, the standard subset of $\Complex^{\times}$
is the union of the interior of the first quadrant
and the positive real axis.  An element of $\Int[\vecti]$
in the standard subset will be referred to as a \textbf{standard Gaussian
integer}, or more simply as a \textbf{standard integer} when the context
is clear.
\end{defn}
Because of the units in $\Int[\vecti]$, each nonzero
ideal $\mathscr{I}$ of $\Int[\vecti]$ has precisely one generator which
is a standard integer.   Henceforth, we refer to generator of $\mathscr{I}$
which is a standard integer as the \textbf{standard
generator} of $\mathscr{I}$.
Unless otherwise stated, whenever we write $\mathscr{I}=(m)$,
to indicate the ideal $\mathscr{I}$ generated by an $m\in\Int[\vecti]$, it will
be understood that $m$ is standard.  Conversely, whenever
we write an ideal $\mathscr{I}$ in the form $(m)$, it
will be understood that $m$ is the standard generator
of $\mathscr{I}$.
Thus, for example, since $(1-\vecti)=\vecti^3(1+\vecti)$ with
$1+\vecti$ standard,
we write $\mathscr{I}=:(1-\vecti)\Int[\vecti]$, defined
as the ideal of Gaussian integers
 divisible
by $1-\vecti$, in the form $\mathscr{I}=(1+\vecti)$.

Similar comments apply to Gaussian primes, factorization, and greatest
common divisor in $\Int[\vecti]$.  By a ``prime in $\Int[\vecti]$", we will
always mean a \textit{standard prime}.  By ``prime factorization"
in $\Int[\vecti]$ we will always mean \textit{factorization
into a product of standard primes}, multiplied by the appropriate unit factor.
Note that the convention regarding standard primes uniquely
determines the unit factor in a prime factorization.  For example,
since
\startdisp
2=\vecti^{3}(1+\vecti)^2
\finishdisp
and $(1+\vecti)^3$ is standard, the above expression is the standard
factorization of the Gaussian integer $2$, and
$\vecti^{3}$ is uniquely determined as the \textit{standard unit factor}
in the prime factorization of $2\in\Int[\vecti]$.

By convention, unless stated otherwise, the ``trivial ideal"
$\Int[\vecti]$ will be understood to belong to the set of ideals
of $\Int[\vecti]$.  The standard generator of
the trivial ideal $\Int[\vecti]$ is, of course, $1$.

To facilitate the statement of Proposition
\ref{prop:inversematrixexplicitdescription}, we
estblish the following conventions.  First, we use $\omega_8$
to denote the unique primitive eighth root of unity in the standard set
of $\Complex^{\times}$.  Observe that
\starteqn\label{eqn:omegaeightdefn}
\omega_8=\frac{\sqrt{2}}{2}(1+\vecti),\quad\text{and}\quad
\omega_8^2=\vecti.
\finisheqn

\vspace*{0.3cm}\noindent\boldmath\textbf{The $\SL{2}{\Int[\vecti]}$-space
$\mathrm{M}_2^{\mathrm{N}}$}.\unboldmath
\begin{defn}  For $N\in\Int[\vecti]$, \boldmath
$\mathrm{M}_2^{\mathrm{N}}\;$\unboldmath
will denote the subset of $\mathrm{Mat}_2(\Int[\vecti])$ consisting
of the elements with determinant $N$.  Since the group $\SL{2}{\Int[\vecti]}$ acts
on $\MtwoN$ by multiplication on the left, $\MtwoN$ is a
$\SL{2}{\Int[\vecti]}$-space.
\end{defn}

It is not difficult to see that the
action of $\SL{2}{\Int[\vecti]}$ on $\MtwoN$ fails to be
transitive, so $\MtwoN$ is not a $\SL{2}{\Int[\vecti]}$-homogeneous
space.  The purpose of the subsequent definitions and results
is to give a description of the orbit structure of the
$\SL{2}{\Int[\vecti]}$-space $\MtwoN$.

Let
\starteqn\label{eqn:omegamdefn}
\Omega_{y}:=\text{a fixed set of representatives of
$\Int[\vecti]/(y)$},\;\text{for all}\;y\in\Int[\vecti].
\finisheqn
It is clear that, for each
$y\in\Int[\vecti]$, there exist a number of possible choices
for $\Omega_{y}$.
For the general result,
Proposition \ref{prop:heckedecomp}, below, the choice
of $\Omega_y$ does not matter, and we leave it unspecified.
However, in the specific applications
of Proposition \ref{prop:heckedecomp}, where $y$ is always
of the form $y=(1+\vecti)^n$ for $n$ a positive integer,
it will be essential
to give an $\Omega_y$ explicitly, which we now do.

So let $n\in\Natural$, $n\geq 1$.
In the definition of $\Omega_{(1+\vecti)^n}$,
we use the ``ceiling" notation, defined by
\startdisp
\lceil q\rceil\;=\text{smallest integer $\geq q$, for $q\in\Rational$.}
\finishdisp
Now set
\starteqn\label{eqn:omegaspecific}
\Omega_{(1+\vecti)^n}=\left\{r+s\vecti\;\text{with $r,\,s\in\Int$,
$0\leq r<2^{\lceil\frac{n}{2}\rceil}$, $0\leq s< 2^{n-
\lceil\frac{n}{2} \rceil}$}\right\}.
\finisheqn
The definition is justified by Lemma \ref{lem:standardresiduerep},
below.

\begin{lem}\label{lem:standardresiduerep}
For $n\geq 1$ an integer, let $\Omega_{(1+\vecti)^n}$
be defined as \eqref{eqn:omegaspecific}.  Then
\startdisp
\Omega_{(1+\vecti)^n}\;\text{is a complete set of representatives of
$\Int[\vecti]/\hspace*{-.4mm}\left((1+\vecti)^n)\right)$}\;\text{for all}\; n.
\finishdisp
\end{lem}

\begin{defn}\label{defn:alphamat}  Let $N\in\Int[\vecti]$ be fixed, and
for each $y\in\Int[\vecti]$ let $\Omega_y$ be as in
\eqref{eqn:omegamdefn}.  Define the matrix $\alpha^\mathrm{N}(m,x)
\in \MtwoN$
as follows,
\starteqn\label{eqn:MofNmxdefn}
\alpha^\mathrm{N}(m,x)=\begin{pmatrix}m&x\\
0&\frac{N}{m}\end{pmatrix},\;\text{for}\;m\in\Int[\vecti],\,
m|N,\, x\in\Omega_{\frac{N}{m}}.
\finisheqn
\end{defn}
It is trivial to verify that $\alpha^\mathrm{N}(m,x)$ as given by
\eqref{eqn:MofNmxdefn} indeed has determinant $N$, \textit{i.e.}
$\alpha^\mathrm{N}(m,x)\in \MtwoN$.  The point of Definition
\ref{defn:alphamat} is given by the following proposition.

\begin{prop}\label{prop:heckedecomp}  For $N\in\Int[\vecti]-\{0\}$,
let $\MtwoN$ be the $\SL{2}{\Int[\vecti]}$-space of matrices with entries
in $\Int[\vecti]$ and determinant $N$.
Define the matrices $\alpha^\mathrm{N}(m,x)$ as in \eqref{eqn:MofNmxdefn}.
Then
\starteqn\label{eqn:heckedecomp}
\mathrm \MtwoN=
\bigcup_{\left\{\stackrel{m\in\Int[\vecti]| \; m|N,}
{\frac{N}{m}\;\text{standard}}
\right\}}
\hspace*{-1.07cm}\cdot\hspace*{1cm}
\bigcup_{x\in\Omega_{\frac{N}{m}}}
\hspace*{-.55cm}\cdot\hspace*{0.55cm}
\SL{2}{\Int[\vecti]}\alpha^\mathrm{N}(m,x),
\finisheqn
and \eqref{eqn:heckedecomp} gives the decomposition of
the $\SL{2}{\Int[\vecti]}$-space $\MtwoN$
into distinct $\SL{2}{\Int[\vecti]}$-orbits.
\end{prop}
We now make some comments concerning
the significance of Proposition \ref{prop:heckedecomp}.
First, a statement equivalent to Proposition \ref{prop:heckedecomp}
is that an arbitrary $\alpha\in M_2^{\rm N}$
has a uniquely determined product decomposition of the form
\starteqn\label{eqn:heckeexplicit}\alpha=
\begin{pmatrix}a&b\\
c&d\end{pmatrix}=\begin{pmatrix}p&q\\
r&s\end{pmatrix}\begin{pmatrix}m&x\\0&\frac{N}{m}\end{pmatrix},\;
\text{with}\, m\in\Go,\; m|N,\; \frac{N}{m}\;\text{standard},\; x\in\Omega_
{\frac{N}{m}},\, pr-qs=1.
\finisheqn
The uniqueness is derived from Proposition \ref{prop:heckedecomp}
as follows.  The
second matrix in the product of \eqref{eqn:heckeexplicit}
is uniquely determined by the matrix de because
of the disjointness of the union in \eqref{eqn:heckedecomp}.
The first matrix in the product appearing in \eqref{eqn:heckeexplicit}
is therefore also uniquely determined.

The second remark is that Proposition \ref{prop:heckedecomp}
may be thought of as the Gaussian-integer
version of the decomposition of elements
of $\mathrm{Mat}_2(\Int)$ of fixed determinant $N$,
sometimes known as the Hecke decomposition.  Occasionally
we refer to \eqref{eqn:heckeexplicit} as the \textit{Gaussian}
Hecke decomposition, to distinguish it from this \textit{classical}
Hecke decomposition in the context of the rational integers.
The proof is the same as that of the classical decomposition
except for some care that has to be taken because of the presence
of additional units in $\Int[\vecti]$.
\newcounter{ranrom}\setcounter{ranrom}{7}
For the classical
Hecke decomposition, see page 110, \S\Roman{ranrom}.4, of \cite{langmodforms},
which is the source of our notation for the Gaussian version.

\vspace*{0.3cm}\noindent\boldmath\textbf{Statement
of the Main Result of \S\ref{sec:lattice}}.  \unboldmath
Let $\Xi$ be an arbitrary subset of $\SL{2}{\Int[\vecti]}$.
Suppose, at first, that $\Xi$ is actually a \textit{subgroup}
of $\SL{2}{\Int[\vecti]}$.
Since $\SL{2}{\Int[\vecti]}\alpha^\mathrm{N}(m,x)$ is an $\SL{2}{\Int[\vecti]}$-space, it is
also a $\Xi$-space.  For general subgroups $\Xi$, however, the action of $\Xi$
on $\SL{2}{\Int[\vecti]}\alpha^{\rm N}(m,x)$
fails to be transitive, \textit{i.e.},
$\SL{2}{\Int[\vecti]}\alpha^\mathrm{N}(m,x)$ is
not a $\Xi$-homogeneous space.  We will now describe
the orbit structure of $\SL{2}{\Int[\vecti]}\alpha^\mathrm{N}(m,x)$ for a specific
subgroup $\Xi$.  In order to make the description of the
subgroup and some related subsets of $\SL{2}{\Int[\vecti]}$
easier, we introduce the epimorphism
\startdisp
\mathrm{red}_{1+\vecti}: \SL{2}{\Int[\vecti]}\rightarrow\SL{2}{\Int[\vecti]
/(1+\vecti)}
\finishdisp
by inducing from the reduction map
\startdisp
\mathrm{red}_{1+\vecti}: \Int[\vecti]\rightarrow\Int[\vecti]/(1+\vecti).
\finishdisp
That is, we ``extend" $\mathrm{red}_{1+\vecti}$ from elements
to matrices by setting
\starteqn\label{eqn:matrixredmapdefn}
\mathrm{red}_{1+\vecti}\left(\begin{pmatrix}p&q\\ r&s\end{pmatrix}\right)=
\begin{pmatrix}\mathrm{red}_{1+\vecti}p&\mathrm{red}_{1+\vecti}q\\
\mathrm{red}_{1+\vecti}r&\mathrm{red}_{1+\vecti}s\end{pmatrix}.
\finisheqn
Since $\Omega_{1+\vecti}=\{0,1\}$, we may identify
$\Int[\vecti]/(1+\vecti)$ with $\{0,1\}$.  Similarly
to the convention with $p,q,r,s\in\Int[\vecti]$,
we use $(\overline{p},\overline{q},\overline{r},\overline{s})$
to denote a quadruple of elements of $\Int[\vecti]/(1+\vecti)$
such that
\startdisp
\overline{p}\hspace*{0.5mm}\overline{s}-\overline{r}\hspace*{0.5mm}\overline{q}=1.
\finishdisp
Here are two elements of $\SL{2}{\Int[\vecti]/(1+\vecti)}$ of particular
interest.
\starteqn\label{eqn:xi12elements}
\overline{I}:=\begin{pmatrix}1&0\\0&1\end{pmatrix},\;
\overline{S}:=\begin{pmatrix}0&1\\1&0
\end{pmatrix}\in\SL{2}{\Int[\vecti]/(1+\vecti)}.
\finisheqn
The notation in \eqref{eqn:xi12elements}
is chosen to remind the reader that $\overline{I}=\red_{1+\vecti}(I)$
and $\overline{S}=\red_{1+\vecti}(S)$, where $I,\, S$ are the standard
generators of $\SL{2}{\Int}$, as in \setcounter{ranrom}{6}
\S\Roman{ranrom}.1
of \cite{jol05}.
Since $\overline{S}^2=\overline{I}$,
it is easy to see that $\{\overline{I},\overline{S}\}$ is a subgroup
of $\SL{2}{\Int[\vecti]/(1+\vecti)}$.
Now define
\starteqn\label{eqn:xi12defn}
\Xi_{12}=\red_{1+\vecti}\inv(\{\overline{I},\, \overline{S}\}).
\finisheqn
Since $\red_{1+\vecti}$ is a morphism, $\Xi_{12}$ is a subgroup
of $\SL{2}{\Int[\vecti]}$.

Also, using the epimorphism $\red_{1+\vecti}$
we define the following subsets of $\SL{2}{\Int[\vecti]}$:
\starteqn
\label{eqn:residuematrices}
\begin{aligned}
\Xi_1&=&\red_{1+\vecti}^{-1}\left(\left\{\begin{pmatrix}0&1\\1&1\end{pmatrix},\,
\begin{pmatrix}1&1\\0&1\end{pmatrix}\right\}\right),\\
\Xi_2&=&\red_{1+\vecti}^{-1}\left(\left\{\begin{pmatrix}1&1\\1&0\end{pmatrix},\,
\begin{pmatrix}1&0\\1&1\end{pmatrix}\right\}\right).
\end{aligned}
\finisheqn
(The subscripts on the $\Xi$ of \eqref{eqn:xi12defn}
and \eqref{eqn:residuematrices}
are chosen in order to
remind the reader of the column in which zeros appear
in the matrices of $\red_{1+\vecti}(\Xi)$.)
Since $\SL{2}{\Int[\vecti]/(1+\vecti)}$
consists of the elements $\overline{I},\overline{S}$ and the four elements
appearing on the right-hand side of \eqref{eqn:residuematrices},
and $\red_{1+\vecti}$
is an epimorphism, we have
\starteqn\label{eqn:sl2xidecomp}
\SL{2}{\Int[\vecti]}=\Xi_1\bigcup\hspace*{-.325cm}\cdot\hspace*{.325cm}\Xi_2
\bigcup\hspace*{-.325cm}\cdot\hspace*{.325cm}\Xi_{12}.
\finisheqn
Unlike $\Xi_{12}$, the subsets $\Xi_1$ and $\Xi_2$ of $\SL{2}{\Int[\vecti]}$
are not subgroups.

All three subsets $\Xi$ in \eqref{eqn:xi12defn}
and \eqref{eqn:residuematrices}
though have a description of the following sort,
which gives some insight into the reason for Sublemma \ref{sublem:finerdecomp},
below.
\starteqn\label{eqn:xisetdescription}
\begin{gathered}
\text{For fixed}\;
\begin{pmatrix}\overline{p}&\overline{q}\end{pmatrix},
\begin{pmatrix}\overline{r}&\overline{s}\end{pmatrix}\in
\left\{
\begin{array}{l}
\begin{pmatrix}1&1\end{pmatrix},\vspace*{0.15cm}\\
\begin{pmatrix}1&0\end{pmatrix},\vspace*{0.15cm}\\
\begin{pmatrix}0&1\end{pmatrix}
\end{array}
\right\}\subset(\SL{2}{\Int[\vecti]/(1+\vecti)})^2,\\
\Xi=\red_{1+\vecti}\inv\left(
\left\{\begin{pmatrix}\overline{p}&\overline{q}\\
\overline{r}&\overline{s}\end{pmatrix},\,
\begin{pmatrix}\overline{r}&\overline{s}\\
\overline{p}&\overline{q}
\end{pmatrix}\right\}\right).
\end{gathered}
\finisheqn
For example, we obtain $\Xi_{12}$ by taking
\startdisp
\begin{pmatrix}\overline{p}&\overline{q}\end{pmatrix}=
\begin{pmatrix}1&0\end{pmatrix}\;\text{and}
\begin{pmatrix}\overline{r}&\overline{s}\end{pmatrix}
=\begin{pmatrix}0&1\end{pmatrix}
\finishdisp
in \eqref{eqn:xisetdescription}.

The reason for introducing the subsets $\Xi$
of \eqref{eqn:residuematrices}
is that they allow us, in Sublemma \ref{sublem:finerdecomp} below
to describe precisely the orbit structure of the $\Xi_{12}$-space $\SL{2}{\Int[\vecti]}
\alpha^\mathrm{N}(m,x)$.
\begin{sublem}\label{sublem:finerdecomp}
Using the notation of
\eqref{eqn:MofNmxdefn} and \eqref{eqn:residuematrices}, we have
\starteqn\label{eqn:finerdecomp}
\SL{2}{\Int[\vecti]}\alpha^\mathrm{N}(m,x)=\bigcup_{\Xi\,=\,\Xi_{1},\,\Xi_{2},\Xi_{12}}
\hspace*{-1cm}\cdot\hspace*{0.7cm}\Xi \alpha^\mathrm{N}(m,x).
\finisheqn
Each of the three sets in the union \eqref{eqn:finerdecomp} is closed
under the action, by left-multiplication, of
$\Xi_{12}$ on $\SL{2}{\Int[\vecti]}\alpha^\mathrm{N}(m,x)$ and equals precisely
one $\Xi_{12}$-orbit in the space
$\SL{2}{\Int[\vecti]}\alpha^\mathrm{N}(m,x)$.
\end{sublem}
\begin{prop}\label{prop:inversematrixexplicitdescription}  Let $\conj$
be the morphism from $\SL{2}{\Complex}$ onto $G$
as in \eqref{eqn:imagematrixconj}.  Let $\Gamma=\SO{3}{\Int[\vecti]}$
be the group of integral points of $G$ in the coordinatization
of $G$ induced by the isomorphism \eqref{eqn:concreteorthequiv}.
Let the subsets $\Xi_1,\, \Xi_2,\,\Xi_{12}$ of $\SL{2}{\Int[\vecti]}$ be as
defined in \eqref{eqn:xi12defn} and \eqref{eqn:residuematrices}.
Let the matrices $\alpha^{\mathrm{N}}(m,x)$ be as in \eqref{eqn:MofNmxdefn}.
Let $\omega_8\in\Complex$ be as in
\eqref{eqn:omegaeightdefn}.
Then we have
\starteqn\label{eqn:inversematrixexplicitdescription}
\conj\inv(\Gamma)=\bigcup_{\delta,=0,1}
\hspace*{-0.5cm}\cdot\hspace*{.5cm}
\left(\frac{1}{\omega_8^{\delta}}
\Xi_{12}\alpha^{\vecti^{\delta}}\hspace{-0.5mm}(\vecti^{\delta},0)\bigcup\hspace*{-0.32cm}
\cdot\hspace*{0.32cm}
\left(\bigcup_{\epsilon=0,1}\hspace*{-0.45cm}\cdot\hspace*{.45cm}
\frac{1}{\omega_8^{\delta}(1+\vecti)}
\Xi_{2}\alpha^{2\vecti^{1+\delta}}\hspace*{-0.7mm}
(\vecti^{1+\delta},\vecti^{\epsilon})
\right)\right).
\finisheqn
\end{prop}
\vspace*{.3cm}

\noindent\textbf{Remarks}\begin{itemize}
\item[(a)]  We use $\Int[\omega_8]$ to denote the ring
generated over $\Int$ by $\omega_8$.  By \eqref{eqn:omegaeightdefn}
we have $\Int[\vecti]\subset\Int[\omega_8]$ and $\Int[\omega_8]=
\Int[\omega_8,\vecti]$.  It follows from Proposition
 \ref{prop:inversematrixexplicitdescription} that
$\conj\inv(\Gamma)\subseteq \SL{2}{\Complex}$ is in fact a subset of
$\SL{2}{\Rational(\omega)}$.  More precisely, of the two
parts of the right-hand side of \eqref{eqn:inversematrixexplicitdescription},
we have
\starteqn\label{eqn:firstsubsetofinvim} \frac{1}{\omega_8^{\delta}}
\Xi_{12}\alpha^{\vecti^{\delta}}
(\vecti^{\delta},0)\subseteq\SL{2}{\Int[\vecti,\omega_8]}
\quad\text{for}\;\delta\in\{0,\,1\}, \finisheqn while
\starteqn\label{eqn:secondsubsetofinvim}
\left(\bigcup_{\epsilon=0,1}\hspace*{-0.45cm}\cdot\hspace*{.45cm}
\frac{1}{\omega_8^{\delta}(1+\vecti)}
\Xi_{2}\alpha^{2\vecti^{1+\delta}}(\vecti^{1+\delta},\vecti^{\epsilon})
\right)\subseteq\mathrm{SL}_2\left
(\Int\left[\vecti,\omega_8,\frac{1}{1+\vecti}\right]\right)
\quad\text{for}\;\delta\in\{0,\,1\} \finisheqn
\item[(b)]  One can easily verify that the set on the left-hand
side of \eqref{eqn:firstsubsetofinvim}
is closed under multiplication, while
the set on the left-hand side of
\eqref{eqn:secondsubsetofinvim} is not.  More precisely,
through a rather lengthy calculation,
not included here, one verifies that
\starteqn\label{eqn:twosubsetproductalt}
\text{for $(x,y)$ a pair of elements of the form of
\eqref{eqn:secondsubsetofinvim},
xy is}\begin{cases}\text{of form \eqref{eqn:secondsubsetofinvim}}\\
\quad\quad\text{or}\\
\text{of form \eqref{eqn:firstsubsetofinvim}.}
\end{cases}
\finisheqn
with each possibility in \eqref{eqn:twosubsetproductalt} being
realized for an appropriate pair $(x,y)$.
These calculations amount to a brute-force
verification of the fact that the right-hand
side of \eqref{eqn:inversematrixexplicitdescription} is closed under multiplication.
But, because $\Gamma$ is a group and $\conj$ a morphism, this fact
also follows from
Proposition \ref{prop:inversematrixexplicitdescription}.
\end{itemize}

The explicit representation
of $\conj\inv(\Gamma)$ in \ref{prop:inversematrixexplicitdescription}
allows us to read off certain group-theoretic facts relating
$\conj\inv(\Gamma)$ to $\SL{2}{\Int[\vecti]}$.  In Lemma \ref{lem:indexlemma}
below we use the notation
\startdisp
[G:H]\;\text{is the index of $H$ in $G$, for any group $G$ with subgroup $H$}.
\finishdisp
\begin{lem}\label{lem:indexlemma}  Let $\conj\inv(\Gamma)$ be the subgroup
of $\SL{2}{\Complex}$ described above,
given explicitly in matrix form
in \eqref{eqn:inversematrixexplicitdescription}.  All the
other notation is also as in
Proposition \ref{prop:inversematrixexplicitdescription}.
\begin{itemize}\item[(a)]  We have
\startdisp
\conj\inv(\Gamma)\cap\SL{2}{\Int[\vecti]}=\Xi_{12}.
\finishdisp
\item[(b)]  We have
\starteqn\label{eqn:index6}
[\conj\inv(\Gamma):\Xi_{12}]=6,\quad [\SL{2}{\Int[\vecti]}:\Xi_{12}]=3.
\finisheqn
Explicitly, the six right cosets of $\Xi_{12}$ in $\conj\inv(\Gamma)$
are the two cosets obtained by letting $\delta$ range
over $\{0,1\}$ in
\startdisp
\frac{1}{\omega_8^{\delta}}
\Xi_{12}\alpha^{\vecti^{\delta}}\hspace{-0.5mm}(\vecti^{\delta},0)=
\frac{1}{\omega_8^{\delta}}\Xi_{12}\begin{pmatrix}
\vecti^{\delta}&0\\0&1\end{pmatrix}
\finishdisp
and the four cosets obtained by
letting $\delta, \epsilon$ range over $\{0,1\}$ independently in
\startdisp
\frac{1}{\omega_8^{\delta}(1+\vecti)}
\Xi_{12}\begin{pmatrix}1&0\\1&1\end{pmatrix}\alpha^{2\vecti^{1+\delta}}\hspace*{-0.7mm}
(\vecti^{1+\delta},\vecti^{\epsilon})=
\frac{1}{\omega_8^{\delta}(1+\vecti)}
\Xi_{12}\begin{pmatrix}\vecti^{1+\delta}&\vecti^{\epsilon}\\
\vecti^{1+\delta}&2+\vecti^{\epsilon}\end{pmatrix}.
\finishdisp
\end{itemize}
\end{lem}

\section{Good Grenier fundamental domains
for arithmetic groups $\Gamma\in\mathrm{Aut}^+(\bbH^3)$}
\label{sec.explicitfd}
\label{sec:explicitfd}
We begin with the following definition, which is fundamental to everything
that follows.
\vspace*{0.3cm}

\noindent\textbf{Definition.}\hspace*{0.2cm}  Let $X$ be a topological space.
Suppose that $\Gamma$  is a group acting topologically on $X$,
\textit{i.e.}, $\Gamma\subseteq\Iso(X)$.
A subset $\scrF$ of $X$ is called an \textbf{exact fundamental
domain for the action of \boldmath $\Gamma$ on $X$\unboldmath} if the following
conditions are satisfied
\begin{itemize}
\item[\textbf{FD 1.}]  The $\Gamma$-translates of $\scrF$ cover $X$, \textit{i.e.},
\startdisp
X=\Gamma \scrF.
\finishdisp
\item[\textbf{FD 2.}]  Distinct $\Gamma$-translates of $\scrF$ intersect only
on their boundaries, \textit{i.e.},
\startdisp
\gamma_1,\gamma_2\in\Gamma,\, \gamma_1\neq\gamma_2\;\text{implies}\;
\gamma_1\scrF\cap
\gamma_2\scrF\subseteq\gamma_1\partial\scrF,\,\gamma_2\partial\scrF.
\finishdisp
\end{itemize}
Henceforth, we will drop the word \textbf{exact} and refer to such an $\scrF$
simply as a \textbf{fundamental domain}.

For the current section, \S\ref{sec:explicitfd}, only,
$G$, instead of denoting $\SO{3}{\Complex}$, will denote $\SL{2}{\Complex}$.
Likewise, instead of denoting $\SO{3}{\Int[\vecti]}$ or
$\conj\inv(\SO{3}{\Int[\vecti]})$,
$\Gamma$ will denote an arbitrary subgroup of $\SL{2}{\Complex}$,
satisfying certain conditions to be given below.  The main
examples to keep in mind are, first, $\Gamma=\SL{2}{\Int}$, the integer
subgroup of $\SL{2}{\Complex}$, and, second,
$\Gamma=\conj\inv(\SO{3}{\Int[\vecti]})$,
the inverse image of the integer subgroup of $\SO{3}{\Complex}$,
described explicitly as a group of fractional
linear transformations in Proposition
\ref{prop:inversematrixexplicitdescription}.

\vspace*{0.3cm}
\noindent\bmth\textbf{Iwasawa decomposition of $\SL{2}{\Complex}$}.\ubmth
\hspace*{0.3cm}For the reader's convenience, we recall only those results
in the context of $\SL{2}{\Complex}$ which we need to proceed.
For proofs and the statements for $\SL{n}{\Complex}$, see
the ``Notation and Terminology" section of \cite{jol06}.  Let
\startdisp
\begin{aligned}
U&=&&\text{upper triangular unipotent matrices in $\SL{2}{\Complex}$, so
$U=\left\{\left.\begin{pmatrix}1&x\\0&1\end{pmatrix}\;\right|\;x\in\Complex
\right\}$},\\
A&=&&\text{diagonal elements of $\SL{2}{\Complex}$ with positive
diagonal entries, so $A=\left\{\left.\begin{pmatrix}y&0\\0&y\inv\end{pmatrix}
\;\right|\;y\in\Real_{+}\right\}$},\\
K&=&&\text{$\mathrm{SU}(2)$, so $K=\{k\in\SL{2}{\Complex}\;|\;kk^*=1\}$}.
\end{aligned}
\finishdisp
Here $x^*$ denotes the conjugate-transpose $\overline{x}^t$ of $x$.

\textit{We have the \textbf{Iwasawa decomposition}}
\startdisp
\SL{2}{\Complex}=UAK,
\finishdisp
\textit{and the product map $U\times A\times K\rightarrow UAK$ is a differential
isomorphism.}

The Iwasawa decomposition induces a system of coordinates $\phi$
on the symmetric space $\SL{2}{\Complex}/K$.  The mapping $\phi$
is a  diffeomorphism
between $\SL{2}{\Complex}/K$ and $\Real^3$.  The details are as follows.
The Iwasawa decomposition gives a uniquely determined
product decomposition of $gK\in\SL{2}{\Complex}/K$ as
\startdisp
gK=u(g)a(g)K,\;\text{where}\, u(g)\in U,\, a(g)\in A\;
\text{are uniquely determined by $gK$}
\finishdisp
Define the \textbf{Iwasawa coordinates} $x_{1}(g)$, $x_{2}(g)\in\Real$,
$y(g)\in\Real^+$ by the relations
\startdisp
u(g)=\begin{pmatrix}1&x_{1}(g)+\vecti x_{2}(g)\\
0&1\end{pmatrix}\,\quad a(g)=\begin{pmatrix}y(g)^{\half}&0\\
0&y(g)^{-\half}\end{pmatrix}.
\finishdisp
By the Iwasawa decomposition, the Iwasawa coordinates of $g$
are uniquely determined.  We emphasize that while $x_{1}(g)$ and $x_{2}(g)$
range over all the real numbers, $y(g)$ ranges over the positive
numbers.   As functions on $G$, $x_1$ $x_2$, and $y$ are invariant
under right-multiplication by $K$.  Thus $x_1$, $x_2$, and $y$
induce coordinates on $G/K$.
Now define the coordinate mappings
$\phi_i:\SL{2}{\Complex}/K\rightarrow \Real$, for $i=1,2,3$, by
\starteqn\label{eqn:phicoordsdefn}
\phi_1=-\log y,\;\phi_2=x_{1},\;\phi_3=x_{2},
\finisheqn
and set
\startdisp
\phi=(\phi_1,\phi_2,\phi_3): G/K\rightarrow\Real^3.
\finishdisp
The mapping $\phi$ is a diffeomorphism of $G/K$ onto $\Real^3$,
because the Iwasawa coordinate system is a diffeomorphism, as is $\log$.
Thus, there exists the inverse diffeomorphism
\startdisp
\phi\inv: \Real^3\rightarrow G/K.
\finishdisp
By \eqref{eqn:phicoordsdefn}, we can write, explicitly,
\starteqn\label{eqn:phiinverseexplicit}
\phi\inv(t_1,t_2,t_3)=t_2+t_3\vecti+e^{-t_1}\vectj,\forallindisp\;
t=(t_1,t_2,t_3)\in\Real^3.
\finisheqn

\vspace*{0.3cm}
\noindent\bmth\textbf{The quaternion model and the coordinate
system on $\SL{2}{\Complex}/K$}.\ubmth \hspace*{0.3cm}
We will use the model $G/K$ as the upper half-space $\bbH^3$,
defined as the following subset of the quaternions.
\starteqn\label{eqn:quatmodeldefn}
\bbH^3=\{x_{1}+x_{2}\vecti+y\vectj,
\;\text{where}\; x_{1}, \,x_{2}\in\Real,\; y\in\Real^+\}.
\finisheqn
Recall that $\SL{2}{\Complex}$ acts transitively on $\bbH^3$ by fractional
linear transformations.  See \S\Roman{ranrom}.0 of \cite{jol05} for
the details
of the action.  We note the relation
\starteqn\label{eqn:fraclintrans}
g\vectj=x_1(g)+x_2(g)\vecti+y(g)\vectj.
\finisheqn
As a result of \eqref{eqn:fraclintrans} and the Iwasawa decomposition,
we may identify $\SL{2}{\Complex}/K$ with $\bbH^3$.  So $\phi:G/K\rightarrow
\Real^3$ induces a diffeomorphism
\startdisp
\phi:\bbH^3\stackrel{\cong}{\longrightarrow}\Real^3.
\finishdisp
Because of \eqref{eqn:fraclintrans}, if $g$ is any element of $G$
such that $g\vectj=z$, then $\phi(g)=\phi(z)$.  Further, beause
of the way we set up the coordinates on $\bbH^3$,
$\phi: \bbH^3\rightarrow\Real^3$ is given explicitly by the same
formulas as \eqref{eqn:phicoordsdefn}.

As explained in, for example, \S\Roman{ranrom}.0 of \cite{jol05},
the kernel of the action of $\SL{2}{\Complex}$ on $\bbH^3$
is precisely the set $\{\pm I\}$, consisting of the identity
matrix and its negative.

For any oriented manifold $X$ equipped with a metric, use the notation
\startdisp
\mathrm{Aut}^+(X)=\;\text{group of
orientation-preserving isometric automorphisms
of $X$.}
\finishdisp
It is a fact that every element of $\mathrm{Aut}^+(\bbH^3)$ is realized
by a fractional linear transformation in $\SL{2}{\Complex}$,
unique up to multiplication by $\pm 1$.  Therefore,
the action of $\SL{2}{\Complex}$ on $\bbH^3$
by fractional linear transformations induces an isomorphism
\starteqn\label{eqn:autometriesiso}
\SL{2}{\Complex}/\{\pm I\}\cong \mathrm{Aut}^+(\bbH^3).
\finisheqn

\vspace*{0.3cm}
\noindent\bmth\textbf{The stabilizer in $\Gamma$
of the first $j$ $\phi$-coordinates.}
\ubmth\hspace*{0.3cm}
In all that follows, if $i,j\in\Natural$, the notation $[i,j]$ is used to
denote the interval of \textit{integers} from $i$ to $j$,
inclusive.  The interval
$[i,j]$ is defined to be the empty set if $i>j$.
\begin{defn} \label{defn:projections} For $i,j\in\{1,2,3\}$, with $i\leq j$,
let \bmth$\phi_{[i,j]}$
be the \textbf{projection of $\bbH^3$ onto the $[i,j]$
factors of $\Real^3$}.  In other words,
we let
\startdisp
\phi_{[i,j]}=(\phi_i,\phi_{i+1},\ldots,\phi_j).
\finishdisp
\end{defn}
Since $\phi$ is a diffeomorphism of $\bbH^3$, $\phi_{[i,j]}$
is an smooth epimorphism of $\bbH^3$ onto $\Real^{i-j+1}$.

If $\scrK$ is any subset of $\{1,2,3\}$, of size $|\scrK|$, then we
can generalize in the obvious way to define the smooth epimorphism
\startdisp
\phi_{\scrK}: \bbH^3\rightarrow\Real^{|\scrK|}.
\finishdisp

Let $\Gamma$ be a group acting by diffeomorphisms
of $\bbH^3$.  For $\gamma\in\Gamma$ we also use $\gamma$ to denote the
diffeomorphism of $\bbH^3$ defined by the left action of $\gamma$
on $\bbH^3$.
Therefore, for $l\in\{1,\ldots 3\}$ the composition
$\phi_l\circ\gamma$ is the $\Real$-valued function on $\bbH^3$ defined by
\startdisp
\phi_l\circ\gamma(z)=\phi_l(\gamma z)\quad\text{for all}\; z\in \bbH^3.
\finishdisp

We use $\Gamma^{\phi_{[1,j]}}$ to denote the subgroup of $\Gamma$
whose action stabilizes the first $i$ coordinates.  In other words,
we set
\startdisp
\Gamma^{\phi_{[1,j]}}=\{\gamma\in\Gamma\;|\; \phi_{[1,j]}
=\phi_{[1,j]}\circ\gamma\}.
\finishdisp
We extend the definition of $\Gamma^{\phi_{[1,j]}}$ to $j=0,4$,
by adopting the conventions
\startdisp
\Gamma^{\phi_{[1,0]}}=\Gamma,\quad\text{and}\quad\Gamma^{\phi_{[1,4]}}=1.
\finishdisp
Note that, by definition, we have the descending sequence of
groups
\startdisp
\Gamma=\Gamma^{\phi_{[1,0]}}\geq
\Gamma^{\phi_1}\geq \Gamma^{\phi_{[1,2]}}\geq
 \Gamma^{\phi_{[1,3]}}\geq\Gamma^{\phi_{[1,4]}}=1.
\finishdisp
Note that the penultimate group
in this sequence, namely $\Gamma^{\phi_{[1,3]}}$,
equals, by definition, the kernel of the action of $\Gamma$ on $\bbH^3$.
Assuming that $\Gamma\subset\SL{2}{\Complex}$, \textit{i.e.} that
$\Gamma$ consists of fractional linear transformations, we always
have
\starteqn\label{eqn:kernelisverysmall}
\Gamma^{\phi_{[1,3]}}=\Gamma\cap\{\pm 1\}.
\finisheqn

Because the $\Gamma^{\phi_{[1,j]}}$ form a descending sequence, for
$k,j\in\{1,2,3\}$ with $k<j$,
we can consider
the left cosets of $\Gamma^{\phi_{[1,k]}}$ in $\Gamma^{\phi_{[1,j]}}$.
The left cosets are the sets of the form
$\Gamma^{\phi_{[1,j]}}\gamma_k$ for $\gamma_k\in\Gamma^{\phi_{[1,k]}}$.
Now let $i,j,k\in\{1,2,3\}$, $l\leq j$, $k<j$.
By the definition of $\Gamma^{\phi_{[1,j]}}$,
the function $\phi_l\circ\gamma_k$ depends only only on the left
$\Gamma^{\phi_{[1,j]}}$-coset to which $\gamma_k$ belongs.
Therefore, for fixed $z$
we may consider $\phi_l\circ\gamma_k(z)$ to be a well-defined function on
the set of left cosets $\Gamma^{\phi_{[1,j]}}\gamma_k$ of
$\Gamma^{\phi_{[1,k]}}$ in $\Gamma^{\phi_{[1,j]}}$.  We may
therefore, speak of the $\Real$-valued function
$\phi_l\circ\Gamma^{\phi_{[1,j]}}\gamma_k$.

In what follows we
will most often apply the immediately preceding paragraph when
$l=j$, and $k=j-1$.  For $\gamma\in\Gamma^{\phi_{[1,j-1]}}$
and $\Delta$ an arbitrary subset of $\Gamma^{\phi_{[1,j]}}$,
we have
\starteqn\label{eqn:Deltaonleftinphij}
\phi_j(\Delta\gamma z)=\{\phi_j(\gamma z)\}.
\finisheqn
therefore, by setting
\startdisp
\phi_j\circ\Gamma^{[1,j]}\gamma(z)=\phi_j(\gamma z),
\finishdisp
we obtain a well-defined function
\startdisp
\phi_j\circ\Gamma^{\phi_{[1,j]}}\gamma:\bbH^3\rightarrow\Real.
\finishdisp
The function $\phi_j\circ\Gamma^{\phi_{[1,j]}}\gamma$
depends only on the $\Gamma^{\phi_{[1,j]}}$-coset to which
$\gamma$ belongs.

For $\gamma\in\Gamma^{\phi_{[1,j-1]}}$, the $\Real$-valued function
$\phi_j\circ\Gamma^{\phi_{[1,j]}}\gamma$
gives the effect of the action of $\gamma\in\Gamma^{\phi_{[1,j-1]}}$
on the $j^{\rm th}$ coordinate of a point.  It
is clear from the definition that
\starteqn\label{eqn:identitycosetdefn}\text{
$\phi_j=\phi_j\circ\gamma$ if and only if $\Gamma^{\phi_{[1,j]}}\gamma$ is the identity left
coset of $\Gamma^{\phi_{[1,j]}}$ in $\Gamma^{\phi_{[1,j-1]}}$}.
\finisheqn

\vspace*{0.3cm}
\noindent\bmth\textbf{Sections of Projections and induced actions
of $\Gamma$.}
\ubmth\hspace*{0.3cm}
As before, suppose that $\Gamma$ is a group acting by diffeomorphisms
on $\bbH^3$, and let $\Gamma^{\phi_{[1,j]}}$ for $j\in\{1,2,3\}$
be defined
as above.

 For any subset $\scrK$ of
the interval of integers $[1,3]$, we let $\scrK^c=[1,3]-\scrK$
be the \textit{complement of $\scrK$ \textit in $[1,3]$}.

\begin{defn} \label{defn:independentof} Let $f$ be a real-valued function
\startdisp
f: \bbH^3\rightarrow \Real.
\finishdisp
Let $\scrK$ a subset of $[1,3]$.  We say that $f$ is \bmth\textbf{independent
of the $\scrK$ coordinates\hspace*{0.15cm}}\ubmth if for every $x,y\in\bbH^3$,
\startdisp
\phi_{\scrK^c}(x)=\phi_{\scrK^c}(y)\;\text{implies}
\;f(x)=f(y).
\finishdisp
\end{defn}

In other words, $f$ is independent of the coordinates in $\scrK$
if and only if $f$ is constant on the fibers of the projection
$\phi_{\scrK^c}$ onto the $\Real$-factors indexed by $\scrK^c$.

For the next observation, we need to introduce the notion
of a section of a projection $\phi_{\scrK}$.  It will not really
matter which section we use, so for simplicity, we choose the zero section.
For a subinterval $[i,j]$ of $\{1,2,3\}$ of size $j-i+1$, define
\startdisp
\sigma_{[i,j]}^0: \Real^{j-1+1}\rightarrow\bbH^3
\finishdisp
by
\startdisp
\sigma_{[i,j]}^0(x_1,\ldots,x_{j-i+1})=(\underbrace{0,\ldots,0}_{i-1},
x_1,\ldots, x_{j-i+1},\underbrace{0,\ldots,0}_{3-j}).
\finishdisp
The map $\sigma_{[i,j]}^0$ is called the \bmth
\textbf{zero section of the projection $\phi_{[i,j]}$}\ubmth.
The terminology comes from the relation
\starteqn\label{eqn:sectionrelation}
\phi_{[i,j]}\sigma_{[i,j]}^0=\Id_{\Real^{i-j+1}},
\finisheqn
which is immediately verified.
The concept of the zero section of the projection can
 be generalized from the case of a projection
associated with an interval $[i,j]$ to that
of an arbitrary subset $\scrK$ of $\{1,2,3\}$, in the obvious way,
although we will not have any use for this generalization in the present
context.

By use of the zero section, we are able to
make a useful reformulation of the condition that $f:\bbH^3\rightarrow\Real$
is independent of the first $j-1$ coordinates.  Let $j\in\{2,3\}$
and $f$ a real values function on $\bbH^3$.  Then
\starteqn\label{eqn:independencereformulation}\text{
$f$ is independent of the first $j-1$ coordinates if and only if
$f\,\sigma^0_{[j,3]}\phi_{[j,3]}=
\sigma^0_{[j,3]}\phi_{[j,3]}\,f$.
}
\finisheqn

The reformulation \eqref{eqn:independencereformulation} allows
us to prove the following result.
\begin{lem}\label{lem:action1cons} Let $\Delta$ be a group
acting on $\bbH^3$, and for $j\in\{1,2,3\}$, let $\phi_{[j,3]}$
be the projection of $\bbH^3$ onto the last $3-j+1$-coordinates
and let $\sigma^0_{[j,3]}$ be the zero section of $\phi_{[j,3]}$.
Suppose that for all $l\in[j,3]$ and $\delta\in\Delta$ the
functions $\phi_l\circ\delta$ are independent of the first
$j-1$ coordinates.  Then $\Delta$ has an induced action on $\Real^{3-j+1}$
defined by \starteqn\label{eqn:sectioninducedact}
\delta_{[j,3]}(\mathbf{t})=\phi_{[j,3]}(\delta\sigma^0_{[j,3]}(\mathbf{t}))
,\;\forallindisp
\;\mathbf{t}=(t_1,\ldots,t_{3-j+1})\in\Real^{3-j+1}. \finisheqn
\end{lem}

It is an immediate consequence of the definitions that
for any group $\tilde{\Gamma}$ acting on $\bbH^3$
by diffeomorphisms, and any subgroup $\Gamma$
of $\tilde{\Gamma}$, we have, for $1\leq i\leq j\leq 3$,
\starteqn\label{eqn:gammastabilizerintersection}
\Gamma^{\phi_{[i,j]}}=(\tilde{\Gamma})^{\phi_{[i,j]}}\cap\Gamma.
\finisheqn
Applying \eqref{eqn:gammastabilizerintersection} to the case
of $\tilde{\Gamma}=\SL{2}{\Complex}$ and $i=1$, we deduce that
\starteqn\label{eqn:stabilizerintersection}
\Gamma^{\phi_{[1,j]}}=\Gamma\cap\SL{2}{\Complex}^{\phi_{[1,j]}},
\finisheqn
for any subgroup $\Gamma\subseteq\SL{2}{\Complex}$.
Because of \eqref{eqn:stabilizerintersection} it is very useful
to have an explicit expression for
$\SL{2}{\Complex}^{\phi_1}$.  We carry out the calculation
using the relations of \eqref{eqn:phicoordsdefn}.

Let $z\in\bbH^3$ with
\startdisp
z=x_1+x_2+y\vectj,
\finishdisp
as in \eqref{eqn:quatmodeldefn}.  Let
\startdisp
g\in\SL{2}{\Complex}\;\text{with}\;
g=\begin{pmatrix}a&b\\c&d\end{pmatrix}.
\finishdisp
Define
\starteqn\label{eqn:ycdzdefn}
y(c,d;z)=
\frac{y(z)}{||c z+d||^2},
\finisheqn
where in \eqref{eqn:ycdzdefn} and from now on,
for a quaternion $z$,
$||z||^2$ denotes the squared norm of a $z$, so that
$||z||^2=z\overline{z}$.
Then we have
\starteqn\label{eqn:ycdzrelation}
y(gz)=y(c,d;z).
\finisheqn
For the details of such calculations, see
\S\Roman{ranrom}.0 of \cite{jol05}.  Since
\startdisp
\phi_1:\bbH^3\rightarrow\Real\;\text{is defined
as}\;-\log y(\cdot),
\finishdisp
and $\log$ is injective, \eqref{eqn:ycdzrelation} implies
that
\starteqn\label{eqn:ginstabilizerprelim}\text{
$g\in\SL{2}{\Complex}^{\phi_1}$ if and only if
$y(c,d;z)=y(z)$ for all $z\in\bbH^3$.}
\finisheqn

By \eqref{eqn:ginstabilizerprelim} and \eqref{eqn:ycdzdefn}, we have
\starteqn\label{eqn:ginstabilizerprelim2}
g\in\SL{2}{\Complex}^{\phi_1}\;\text{if and only if}\;
||cz+d||^2=1\forallindisp z\in\bbH^3.
\finisheqn
Clearly, the condition $||cz+d||^2=1$ is satisfied
for all $z\in\bbH^3$ if and only if $c=0$ and $||d||=1$.  We therefore
deduce from \eqref{eqn:ginstabilizerprelim2} that
\starteqn\label{eqn:firstcoordfixing}
\SL{2}{\Complex}^{\phi_{1}}=\left\{\left.
\begin{pmatrix}\omega\inv&x\\0&\omega\end{pmatrix}\;\right|\;
x,\,\omega\in\Complex,
\,||\omega||=1\right\}.
\finisheqn

As a result of \eqref{eqn:firstcoordfixing}, we can
easily verify that for $\gamma\in\Gamma^{\phi_1}$, $l\in[2,3]$,
the functions $\phi_l\circ\delta$ are independent
of the first coordinate.  So we can apply Lemma \ref{lem:action1cons},
in this case, with $j=2$ and deduce that
\begin{lem}\label{lem:action1consapptogamma}
Let $\Gamma\subseteq\mathrm{Aut}^+(\bbH^3)$,
$\phi_{[2,3]}$ be the projection of $\bbH^3$ onto the last $2$ coordinates
and let $\sigma^0_{[2,3]}$ be the zero section of $\phi_{[2,3]}$.
Then $\Gamma$ has an induced action on $\Real^{2}$
defined by \starteqn\label{eqn:sectioninducedactapplied}
\gamma_{[2,3]}(\mathbf{t})=\phi_{[2,3]}(\gamma\sigma^0_{[2,3]}(\mathbf{t}))
,\;\forallindisp
\;\mathbf{t}=(t_1,t_2)\in\Real^{2}. \finisheqn
\end{lem}

The following Theorem is a special case of the main result
of the first chapter of \cite{brennerthesis}.
\begin{thm}\label{thm:gammagoodgrenier}
Let $\Gamma$ be a subgroup of $\SL{2}{\Complex}$, acting on $\bbH^3$
on the left by fractional linear transformations.
Suppose that $\Gamma$ is commensurable to $\SL{2}{\Int[\vecti]}$.
 Let $\scrG$ be a fundamental domain for the induced action of
$\Gamma^{\phi_{[2,3]}}/\{\pm 1\}$ on $\Real^2$.  Assume further
that $\scrG=\overline{\intrr(\scrG)}$.  Define
\starteqn\label{eqn:scrF1gammadefn}
\scrF_1=\{z\in\bbH^3\;|\; \phi_1(z)\leq \phi_1(\gamma z),\;
\text{for all}\; \gamma\in\Gamma\}.
\finisheqn
Set
\starteqn\label{eqn:scrFdefn}
\scrF(\scrG)=\phi\inv_{{[2,3]}}(\scrG)\cap\scrF_1.
\finisheqn
\begin{itemize}
\item[(a)]  We have $\scrF(\scrG)$ a fundamental domain for the action
of $\Gamma/\{\pm 1\}$ on $\bbH^3$.
\item[(b)]  We have
\starteqn\label{eqn:scrFtopology}
\scrF(\scrG)=\overline{\intrr\big(\scrF(\scrG)\big)}.
\finisheqn
\item[(c)]
Further, $\intrr\big(\scrF_1\big)$ and $\intrr\big(\scrF(\scrG)\big)$
have explicit descriptions as follows.
\starteqn\label{eqn:scrF1intrrdescription}
\intrr\big(\scrF_1\big)=
\{z\in\bbH^3\;|\; \phi_1(z)< \phi_1(\gamma z),\;
\text{for all}\; \gamma\in\Gamma-\Gamma^{\phi_1}\},
\finisheqn
and
\starteqn\label{eqn:scrFintrr}
\intrr\big(\scrF(\scrG)\big)=\phi_{[2,3]}\inv\big(\intrr(\scrG)\big)\cap
\intrr(\scrF_1),
\finisheqn
\end{itemize}
\end{thm}
Considering the coordinate system $\phi$ on $\bbH^3$ as fixed,
we may think of the fundamental domain $\scrF$ for
$\Gamma^{\phi_{[1,3]}}\backslash\Gamma$ to be a function of the fundamental
domain $\scrG$ for the induced action of $\Gamma^{\phi_1}$
on $\Real^2$.  When we wish to stress this dependence
of $\scrF$ on $\scrG$, we will write $\scrF(\scrG)$ instead of $\scrF$.

\begin{defn}  Suppose that
$\Gamma\subseteq\mathrm{Aut}^+(\bbH^3)$ is commensurable to
$\SL{2}{\Int[\vecti]}$.
Let $\scrG$ be a fundamental domain for the induced
action of $\Gamma^{\phi_{[1,3]}}\backslash\Gamma^{\phi_1}$
on $\Real^2$ satisfying $\scrG=\overline{\intrr(\scrG)}$.
Then the fundamental domain $\scrF(\scrG)$
for the action of $\Gamma^{\phi_{[1,3]}}\backslash\Gamma$
defined in \eqref{eqn:scrFdefn} is called the \bmth\textbf{good
Grenier fundamental domain for the action of $\Gamma$ on $\bbH^3$
associated to the fundamental domain $\scrG$}\ubmth.
\end{defn}
The reference to the fundamental domain $\scrG$ is often omitted in
practice.

Henceforth, we drop the explicit reference to $\Gamma^{\phi_{[1,3]}}$
and speak of a
\textit{fundamental domain of $\Gamma^{\phi_{[1,3]}}\backslash\Gamma$}
as a \textit{fundamental domain of $\Gamma$}.
By \eqref{eqn:kernelisverysmall},
$\Gamma$ is at worst a two-fold cover of
$\Gamma^{\phi_{[1,3]}}\backslash\Gamma$,
so this involves only a minor abuse of terminology.

We will give an expression for a good Grenier fundamental
domain $\scrF(\scrG)$ for $\conj\inv(\SO{3}{\Int[\vecti]})$
in terms of explicit inequalities,
in \eqref{eqn:scrFscrGfirstform}, and again as a
convex polytope in $\bbH^3$, in Proposition \ref{prop:gammafunddom},
below.
\vspace*{0.3cm}

\noindent\bmth\textbf{Example: The Picard domain $\scrF$
for $\SL{2}{\Int[\vecti]}$.}
\ubmth\hspace*{0.3cm}
Define the following rectangle in $\Real^2$:
\starteqn\label{eqn:picardGdescription}
\scrG_{\SL{2}{\Int[\vecti]}^{\phi_1}}
=\left\{(t_1,t_2)\in\Real^2\;\left|\; t_1\in\left[-\half,\half\right],\,
t_2\in\left[0,\half\right]\right.\right\}.
\finisheqn
It is easy to verify, from an explicit description of
$\SL{2}{\Int[\vecti]}^{\phi_1}$, deduced from \eqref{eqn:firstcoordfixing}
that $\scrG_{\SL{2}{\Int[\vecti]}^{\phi_1}}$ is a fundamental
domain for the action of $\SL{2}{\Int[\vecti]^{\phi_1}}/\{\pm 1\}$.

Further, it is obvious that
\startdisp
\scrG_{\SL{2}{\Int[\vecti]}^{\phi_1}}=
\overline{\intrr(\scrG_{\SL{2}{\Int[\vecti]}^{\phi_1}})}.
\finishdisp

Therefore, Theorem \ref{thm:gammagoodgrenier}
applies.  We deduce that, with $\scrF_1$,
$\scrF(\scrG_{\SL{2}{\Int[\vecti]}^{\phi_1}})$ defined as in Theorem
\ref{thm:gammagoodgrenier}, we have
\startdisp\text{
$\scrF:=\scrF(\scrG_{\SL{2}{\Int[\vecti]}^{\phi_1}})$ is a good
Grenier fundamental domain for $\SL{2}{\Int[\vecti]}$.
}
\finishdisp
The fundamental domain $\scrF$ is defined
in \S\Roman{ranrom}.1 of \cite{jol05}, where, in keeping with
classical terminology, $\scrF$ is called the \textbf{Picard domain}.

In order to complete the example, we now give an
explicit description of the set $\scrF_1$, which will
allow the reader to see that ``our" $\scrF$ is exactly
the same as the Picard domain.  It can be shown that
$\scrF_1$ is the subset of $\Real^3$ whose image under
the diffeomorphism $\phi\inv$ is given as follows.
\starteqn\label{eqn:scrF1sl2zidescript}
\phi\inv(\scrF_1)=\{z\in\bbH^3\;|\; ||z-m||\geq 1,\;\text{for all}\,
m\in\Int[\vecti]\}.
\finisheqn

Of the infinite set of
inequalities defining $\scrF_1$, all except the one
with $d=0$, \textit{i.e.} $||z||^2\geq 1$, are
trivially satisfied on
$\phi_{[2,3]}\inv\left(\scrG_{\SL{2}{\Int[\vecti]}^{\phi_1}}\right)$.
Thus, from \eqref{eqn:scrF1sl2zidescript}, \eqref{eqn:picardGdescription},
and \eqref{eqn:scrFdefn}, we recover the description of the Picard domain by
finitely many inequalities given in \S\Roman{ranrom}.1 of \cite{jol05}.
\starteqn\label{eqn:picardexplicitfinite}
\scrF(\scrG_{\SL{2}{\Int[\vecti]}^{\phi_1}})=
\left\{z\in\bbH^3\;\left|\; \, x_1\in\left[-\half,\half\right],\,
x_2\in\left[0,\half\right]
y,\, ||z||^2\geq 1\right.\right\}.
\finisheqn

\section{Explicit description of the fundamental domain
for the action of $\SO{3}{\Int[\vecti]}$ on $\bbH^3$}
We now proceed to consider the special case of
$\conj\inv(\SO{3}{\Int[\vecti]})$ in Theorem
\ref{thm:gammagoodgrenier} above.
In keeping with the general practice of the present paper,
we will go back to using $G$ to denote $\SO{3}{\Complex}$
exclusively, and $\Gamma$ to denote the group
$\SO{3}{\Int[\vecti]}$.  Since we are always
in this section in the setting of subgroups of $\SL{2}{\Complex}$,
we will abuse notation slightly and use $\Gamma$ to denote
the isomorphic inverse image $\conj\inv(\Gamma)$ of
$\Gamma=\SO{3}{\Int[\vecti]}$ in $\SL{2}{\Complex}$.

Also, we treat $\Real^2$, the image of the projection
$\phi_{[2,3]}$, as $\Complex$, by identifying the point $(t_1,t_2)\in\Real^2$
with $t_1+\vecti t_2$.  Thus, our ``new" $\phi_{[2,3]}$ is defined
in terms of the ``old" $\phi$-coordinates by
\starteqn\label{eqn:phi23new}
\phi_{[2,3]}(z)=\phi_2(z)+\vecti\phi_3(z).
\finisheqn

\begin{prop}\label{prop:scrF1forgamma}\bmth
\textbf{First form of $\scrF_1$.}\ubmth \hspace*{1mm} Let $\scrF_1$ be as
defined in \eqref{eqn:scrF1gammadefn}.  All other
notation has the same meaning as in Theorem \ref{thm:gammagoodgrenier}.
Then we have
\starteqn\label{eqn:scrF1gammaexplicit}
\scrF_1=\{z=x(z)+y(z)\vectj\in\bbH^3\;
|\; ||x(z)-d||^2+y(z)^2\geq 2,\;\text{for}\;
d\in1+(1+\vecti)\Int[\vecti]\},
\finisheqn
and $\intrr(\scrF_1)$ is the same as in \eqref{eqn:scrF1gammaexplicit},
but with strict inequality instead of nonstrict inequality.
\end{prop}

\noindent\bmth\textbf{Fundamental domain $\scrG$
for $\Gamma^{\phi_1}$.}\ubmth\hspace*{3mm}
In order to complete the explicit determination of a good
Grenier fundamental domain $\scrF$ for $\Gamma$, it remains
to give describe a suitable fundamental domain $\scrG$ for $\Gamma^{\phi_1}$.
Using \eqref{eqn:stabilizerintersection}, \eqref{eqn:firstcoordfixing},
and the description of $\Gamma$ in
\eqref{eqn:inversematrixexplicitdescription}
we deduce that
\starteqn\label{eqn:gammaphi1explicit}
\Gamma^{\phi_1}=\left\{\left.
\begin{pmatrix}\omega_8^{\delta}&\omega_8^{\delta}b\\
0&\omega_8^{-\delta}\end{pmatrix}\;\right|\; b\in(1+\vecti)\Int[\vecti],\;
\delta\in\{0,1\}\right\}.
\finisheqn
It follows from \eqref{eqn:gammaphi1explicit} that the subgroup
of unipotent elements of $\Gamma^{\phi_1}$ is
\starteqn\label{eqn:gammaphi1Uexplicit}
(\Gamma^{\phi_1})_U=\begin{pmatrix}1&(1+\vecti)\Int[\vecti]\\0&1\end{pmatrix}.
\finisheqn
We make note of certain group-theoretic properties of $\Gamma^{\phi_1}$
and $(\Gamma^{\phi_1})_U$ that are used in determining
the fundamental domains.
First, we define the following generating elements:
\starteqn\label{eqn:gammagensdefn}
R_{\frac{\pi}{2}}=\begin{pmatrix}\omega_8&0\\0&
\omega_8\inv\end{pmatrix},\;
T_{1+\vecti}=\begin{pmatrix}1&1+\vecti\\0&1\end{pmatrix},\;\text{and}\;
T_{1-\vecti}=\begin{pmatrix}1&1-\vecti\\0&1\end{pmatrix}.
\finisheqn
It is easily verified, using \eqref{eqn:gammaphi1explicit}
and \eqref{eqn:gammaphi1Uexplicit}, that
\starteqn\label{eqn:generatingsets}
(\Gamma^{\phi_1})_U=\langle T_{1+\vecti},\, T_{1-\vecti}\rangle,\quad
\Gamma^{\phi_1}=\langle R_{\frac{\pi}{2}}, \,T_{1+\vecti},
\,T_{1-\vecti}\rangle.
\finisheqn
We calculate, from the definition of $R_{\frac{\pi}{2}}$
and \eqref{eqn:generatingsets}, that
\startdisp
\conj(R_{\frac{\pi}{2}})(\Gamma^{\phi_1})_U=(\Gamma^{\phi_1})_U.
\finishdisp
Since $\Gamma^{\phi_1}$ is generated by $(\Gamma^{\phi_1})_U$
and $R_{\frac{\pi}{2}}$, and $R_{\frac{\pi}{2}}$ has order
$4$, we deduce that
\starteqn\label{eqn:Unormality}\text{
$(\Gamma^{\phi_1})_U$ is normal in $\Gamma^{\phi_1}$ with
$[\Gamma^{\phi_1}:(\Gamma^{\phi_1})_U]=4$.}
\finisheqn
Let $T$ be any element of $(\Gamma^{\phi_1})_U$.
Then we have a more precise version of \eqref{eqn:Unormality},
\starteqn\label{eqn:cyclicgpcosetreps}\text{The
group $\langle TR_{\frac{\pi}{2}}\rangle$ of order $4$
is a set of representatives for the coset group $\Gamma^{\phi_1}/(\Gamma^{\phi_1})_U$.}
\finisheqn
Applying \eqref{eqn:cyclicgpcosetreps} to the case $T=T_{1-\vecti}$, we have
\starteqn\label{eqn:cyclicgpcosetrepsT1plusi}\text{The
group $\langle T_{1-\vecti}R_{\frac{\pi}{2}}\rangle$ of order $4$
is a set of representatives for the coset group $\Gamma^{\phi_1}/(\Gamma^{\phi_1})_U$.}
\finisheqn

It is easily verified that the action of $R_{\frac{\pi}{2}}$
on $\Complex$ is rotation by an angle $\pi/2$ about the
fixed point $0$.  Furthermore, calculate from \eqref{eqn:gammagensdefn}
that
\startdisp
T_{1-\vecti}R_{\frac{\pi}{2}}=\conj(T_{1})R_{\frac{\pi}{2}}.
\finishdisp
Therefore,
\starteqn\label{eqn:rotationdescription}\text{The action
of $T_{1-\vecti}R_{\frac{\pi}{2}}$ on $\Complex$
is rotation by $\pi/2$ about $1$.}
\finisheqn
The following
statement is a special case of Lemma 2.2.7 of \cite{brennso3preprint}.
\begin{lem}\label{lem:finiteindexfunddom}
Let $\scrG_U$ be a
fundamental domain for the action of $(\Gamma^{\phi_1})_{U}$
on $\bbH^3$, satisfying
\startdisp
T_{1+\vecti}R_{\frac{\pi}{2}}(\scrG_U)=\scrG_U.
\finishdisp
Let $\scrG$ be a fundamental domain
for the action of $\langle T_{1+\vecti}R_{\frac{\pi}{2}}\rangle$
on $\scrG$.
Then $\scrG$ a fundamental domain for the action of
$\Gamma^{\phi_1}$ on $\bbH^3$.
\end{lem}

In order to define and work with the sets $\scrG_U$ and $\scrG$
which will be fundamental domains for the action of
$\Gamma^{\phi_1}_U$ and $\Gamma^{\phi_1}$, it is useful to introduce
the notion of a convex hull in a totally geodesic metric space.

A metric space $(X,d)$ will be called \textbf{totally geodesic} if
for every pair of points $p_1,p_2\in X$, $p_1\neq p_2$ there is a unique
geodesic segment connecting $p_1,p_2$.  In this situation,
the (closed) geodesic segment connecting
$p_1,p_2$ will be denoted $[p_1,p_2]_{d}$.  A point $x\in X$
is said to lie \bmth\textbf{between $p_1$ and $p_2$}\ubmth\hspace*{0.5mm}
when $x$ lies on $[p_1,p_2]_{d}$.  We then say that
$\scrS\subset X$ is \textbf{convex} when $p_1,\,p_2\in\scrC$
and $p_3$ between $p_1$ and $p_2$ implies
that $p_3\in\scrS$.  Let $p_1,\ldots,p_r$ be $r$ points
in $X$.  The points determine a set
\startdisp
\scrC_{d}(p_1,\ldots, p_r)
\finishdisp
called the \textbf{convex closure} of $p_1,\ldots,p_r$,
described as the smallest convex subset of $X$ containing the set
$\{p_1,\ldots,p_r\}$.

Obviously, we can apply the notion of convex hull to any set $\scrS$,
rather than a finite set of points.  The definition remains the same,
namely that $\scrC_d(\scrS)$ is the smallest convex subset of $X$
containing $\scrS$.  In general we will use the notation
\startdisp
\scrC_d(\scrS_1,\ldots,\scrS_r)=\scrC_d\left(\bigcup_{i=1,\ldots r}
\scrS_i\right).
\finishdisp

In particular, if we apply these notions to $X=\Real^2$
with the ordinary Euclidean metric $\mathrm{Euc}$, then
the geodesic segment $[p_1,p_2]_{\mathrm{Euc}}$ is
just the line-segment joining $p_1,p_2$.  Further,
provided that not all the $p_i$ are collinear, $\scrC(p_1,\ldots p_r)$
is a closed convex polygon whose vertices are located at a subset of
$\{p_1,\ldots,p_r\}$.

We first use the notion of convex closure to record an
elementary facts concerning the fundamental domains
of groups of translations acting on $\Real^2$, identified with $\Complex$
in the usual way.  Let $\omega_1,\omega_2\in\Complex$
be linearly independent over $\Real$.  Then
$\Int\omega_1+\Int\omega_2$ is a lattice in $\Complex$,
and it is well known that all lattices in $\Complex$
are of this form for suitable $\omega_1, \omega_2$
Let $T$ denote the group of translations by elements
of $\Int\omega_1+\Int\omega_2$ acting on $\Complex$.  Then we have
\starteqn\label{eqn:latticetransF}
\scrC(0,\omega_1,\omega_2,\omega_1+\omega_2)\;\text{is a fundamental
domain for the action of $\Int\omega_1+\Int\omega_2$ on $\Complex$}.
\finisheqn
Now we define the following polygons in $\Complex\cong\Real^2$.  Let
\startdisp
\scrG_U=\scrC_{\mathrm{Euc}}(0,2,1+\vecti,1-\vecti),
\finishdisp
and let
\starteqn\label{eqn:scrGdefn}
\scrG=\scrC_{\mathrm{Euc}}(1,2,1+\vecti).
\finisheqn
The relation between the polygons is that $\scrG_U$ is a square
centered at $1$, while $\scrG$ is an isosceles right triangle
inside $\scrG_U$,
with vertices at the center of $\scrG_U$ and two of the corners
of $\scrG_U$.  Therefore, it follows from \eqref{eqn:rotationdescription}
that we have
\starteqn\label{eqn:scrGfunddominscrGU}
\scrG_U=\bigcup_{i=0,1,2,3}(T_{1+\vecti}R_{\frac{\pi}{2}})^i\scrG,\;\text{with}
\;(T_{1+\vecti}R_{\frac{\pi}{2}})^i\scrG\cap\scrG
\subseteq\partial\scrG,\;\text{for}\; i\not\equiv0\mod 4.
\finisheqn
The relations \eqref{eqn:latticetransF} and \eqref{eqn:scrGfunddominscrGU}
lead to the following lemma.

\begin{lem}\label{lem:scrGfunddom}  Let $\Gamma^{\phi_1}$
be as given in \eqref{eqn:gammaphi1explicit} and
$(\Gamma^{\phi_1})_U$ as given in
\eqref{eqn:gammaphi1Uexplicit}.
\begin{itemize}
\item[(a)]  The set $\scrG_U$ is a fundamental domain for the
induced action of $(\Gamma^{\phi_1})_U$ on $\Complex\cong\Real^2$.
\item[(b)]  $\scrG$ is a fundamental domain for the induced
action of $\langle T_{1+\vecti}R_{\frac{\pi}{2}}\rangle$
on $\scrG_U$.
\item[(c)]  The set $\scrG$ is a fundamental domain for the
induced action of $\Gamma^{\phi_1}$ on $\Complex\cong\Real^2$.
\end{itemize}
\end{lem}

\bmth
\noindent\textbf{Form of $\scrF$ in terms of explicit inequalities.}
\ubmth\hspace*{3mm}
Combining Part (c) of Lemma \ref{lem:scrGfunddom},
Proposition \ref{prop:scrF1forgamma}, and \eqref{eqn:scrFdefn},
we deduce that
\startdisp
\scrF(\scrG)=\{z\in\bbH^3\;|\; \phi_{[2,3]}(z)\in\scrC_{\mathrm{Euc}}(1,2,1+\vecti),\;
||x(z)-m||^2+y(z)^2\geq 2,\text{for}\; m\in1+(1+\vecti)\Int[\vecti]\}.
\finishdisp
By \eqref{eqn:phi23new}, the first
condition in the description of $\scrF(\scrG)$ above may be replaced
by
\starteqn\label{eqn:trianglecondition}
x(z)\in\scrC_{\mathrm{Euc}}(1,2,1+\vecti)
\finisheqn
Let $z\in\Complex$ satisfying \eqref{eqn:trianglecondition}.
The element $m=1$ is the element
 of $1+(1+\vecti)\Int[\vecti]$ closest to $x(z)$.
Therefore, for $z$ satisfying
\eqref{eqn:trianglecondition},
the condition
\startdisp
||x(z)-m||^2+y(z)^2\geq 2,\text{for all}\; m\in1+(1+\vecti)\Int[\vecti]
\finishdisp
reduces to $||x(z)-1||^2+y(z)^2\geq 2$.  So we may rewrite
the description of $\scrF(\scrG)$ in the form
\starteqn\label{eqn:scrFscrGfirstform}
\scrF(\scrG)=\{z\in\bbH^3\; x(z)\in\scrC_{\mathrm{Euc}}(1,2,1+\vecti),\;
||x(z)-1||^2+y(z)^2\geq 2\}.
\finisheqn
\vspace*{3mm}

\noindent\bmth
\textbf{Additional facts regarding convex hulls and totally geodesic
hypersurfaces in $\overline{\bbH^3}$.}
\ubmth\hspace*{3mm}
We now extend our ``geodesic hull" treatment of $\scrF$ from the boundary
into the interior of $\bbH^3$.
We first recall certain additional
facts regarding convex hulls and totally geodesic hypersurfaces
in $\bbH^3$.

The description of the geodesics in $\bbH^2$ is well known, but
the corresponding description of the totally geodesic surfaces in
$\bbH^3$ perhaps not as well known, so we recall it here.
Henceforth we abbreviate
``totally geodesic'' by t.g.  Although all t.g. surfaces are
related by isometries, in our model they have two basic types. The
first type is a vertical upper half-plane passing through the
origin with angle $\theta$ measured counterclockwise from the real
axis, which we denote by $\bbH^2(\theta)$.  The second type is an
upper hemisphere centered at the origin with radius $r$, which we
will denote by $\bbS^+_r(0)$.  The t.g. surfaces of $\bbH^3$ are
the $\bbH^2(\theta)$, the $\bbS^+_r(0)$, and their translates by
elements of $\Complex$.  For each of the basic t.g. surfaces, we
produce an isometry $g\in\text{Aut}(\bbH^3)$, necessarily
orientation-reversing, such that $\Fix(g)$ is precisely
the surface in question.  The existence of such a
$g$ shows that the surface is a t.g.
surface.

We define
\startdisp
\overline{\bbH^3}=\bbH^3\cup\Complex\cup\infty
\finishdisp
to be the usual closure of $\bbH^3$ and extend the action of fractional
linear transformations and the notion of the convex hull in
the usual way.  For any subset $\scrS$ of $\bbH^3$, $\overline{\scrS}$
will denote the closure in $\overline{\bbH^3}$.  For $g\in\mathrm{Aut}
(\bbH^3)$, we will likewise use $g$ to denote the extension of $g$
to the closure $\overline{\bbH^3}$.
Henceforth, we will work exclusively in the setting of the closure
$\overline{\bbH^3}$ of $\bbH^3$.  Thus, we will
actually identify the closures of the t.g. surfaces.

The basic orientation-reversing isometry of $\overline{\bbH^3}$
may be denoted $R^*$.
With $x_1+x_2\vecti+y\vectj \in\overline{\bbH^3}$, we have
\startdisp R^*(x_1+x_2\vecti+y\vectj)=x_1-x_2\vecti+y\vectj.
\finishdisp Clearly, we have $\Fix(R^*)=\overline{\bbH^2(0)}$.  To obtain
isometries corresponding to the other vertical planes, let
\startdisp
R_{\theta}=\begin{pmatrix}e^{\vecti\theta/2}&0\\
0&e^{-\vecti\theta/2}\end{pmatrix}. \finishdisp Because
$R_{\theta}\overline{\bbH^2(0)}=\overline{\bbH^2(\theta)}$,
we have \startdisp
\Fix(\conj(R_{\theta})R^*)=\overline{\bbH^2(\theta)}.
\finishdisp

To define the isometry $I$ such that $\Fix(I)$ is the basic
hemisphere $\overline{\bbS_0^+(1)}$, let $\overline{z}$ denote the conjugate
of the quaternion $z$, \textit{i.e.} if $z=x_1+x_2\vecti+y\vectj$
then $\overline{z}=x_1-x_2\vecti-y\vectj$.  For $z\in\overline{\bbH^3}$, set
\startdisp I(z)=1/\overline{z}. \finishdisp We have the equality
$z/I(z)=||z||^2$.  Observe that $\overline{\bbS^+_1(0)}$ is precisely the set
of quaternions in $\overline{\bbH^3}$ of norm one.   Thus,
 $\Fix(I)=\overline{\bbS^+_1(0)}$. For the more general
hemispheres $\overline{\bbS^+_r(0)}$, set \startdisp A(r)=\begin{pmatrix}
\sqrt{r}&0\\
0&\frac{1}{\sqrt{r}}
\end{pmatrix}.
\finishdisp Then, since $A(r)\overline{\bbS^+_1(0)}=
\overline{\bbS^+_r(0)}$, we have
$\Fix(\conj(A(r))I)=\overline{\bbS_r^+(0)}$.

In order to denote the convex hull in $\overline{\bbH^3}$, we use
the notation $\coh$.  Therefore, if $\intd s^2$ is the hyperbolic
metric on $\overline{\bbH^3}$, we have
\startdisp
\coh(p_1,\ldots, p_r)=\mathscr{C}_{\intd s^2}(p_1,\ldots, p_r),
\finishdisp
in terms of our original notational conventions.

Let $p_1,\ldots p_r\in\overline{\bbH^3}$,
for $r>3$ not lying on the same totally
geodesic surface, such that, for each $i$, $1\leq i\leq r$,
\[
p_i\notin\coh(p_1,\ldots,p_{i-1},p_{i+1},\ldots,p_r).
\]
Then the set $\coh(p_1,\ldots, p_r)$ will be called the
\textbf{solid convex polytope with vertices at \boldmath$p_1,\ldots,
p_r$\unboldmath}.  It is clear that for any $p_1,\ldots
p_r\in\overline{\bbH^3}$ not lying in the same totally geodesic
surface, $\coh(p_1,\ldots, p_r)$ is a solid convex polytope with vertices
consisting of some subset of the $r$ points.
\vspace*{3mm}

\noindent\bmth\textbf{Description of $\scrF(\scrG)$ as a solid convex
polytope.}\hspace*{3mm}
\ubmth
\begin{prop}\label{prop:gammafunddom}
The solid convex polytope with four vertices in $\overline{\bbH^3}$
given by
\starteqn\label{eqn:scrFscrGsecondform}
\scrF(\scrG)=\scrC_{\mathbf{H}}(1+\sqrt{2}\vectj,2+\vectj,1+\vecti+\vectj,\infty)
\finisheqn
is a good Grenier fundamental domain for the action of
$\Gamma=\conj\inv(\SO{3}{\Int[\vecti]})$ on $\overline{\bbH^3}$.
\end{prop}

\section{$\SOM{2}{1}_{\Int}$ as a group of fractional
linear transformations}
We will now use the results of \S\ref{sec:lattice} and \S\ref{sec:explicitfd}
to deduce a realization of $\Gamma_{\Int}=\SOM{2}{1}_{\Int}$ as a group
of fractional linear transformations, as well as a description
of a fundamental domain for $\Gamma_{\Int}$ acting on $\bbH^2$
that is in some sense (to be explained precisely below) compatible with the
fundamental domain of $\Gamma$ acting on $\bbH^3$.

We maintain the notational conventions established in
\S\ref{sec:lattice}.  In particular, $G=\SO{3}{\Complex}$
and $\Gamma=\SO{3}{\Int[\vecti]}$.  It is crucial, for the moment,
that we observe the distinction between $G,\Gamma$ and their
isomorphic images under $\conj\inv$.
\begin{defn} \label{defn:gammaz} Set
\starteqn\label{eqn:gammazdefn}
\Gamma_{\Int}=\conj(\SL{2}{\Real}\cap\conj\inv(\Gamma)).
\finisheqn
\end{defn}
\begin{rem}
Note that the elements of $\Gamma_{\Int}$ do not have real entries!
The na\"{i}ve approach to the definition of $\Gamma_{\Int}$
would be to take
the elements of $\Gamma$ with real entries, as in the case of
$\SL{2}{\Int[\vecti]}$ and $\SL{2}{\Int}$.  However, this clearly
cannot be the right definition
because the resulting discrete group would be contained
in $\SOR{3}$, hence compact, and hence finite.
The justification for Definition \ref{defn:gammaz}
is contained in Proposition \ref{prop:gammaziso}, below.
\end{rem}

Recall the orthonormal basis $\beta$ for $\mathrm{Lie}(\SL{2}{\Complex})$
defined at $\eqref{eqn:beta}$.  Define a new basis
$\eta$ by specifying the change-of-basis matrix
\starteqn\label{eqn:etadefn}
\alpha^{\beta\mapsto\eta}=\mathrm{diag}(1,-\vecti,1).
\finisheqn

Let $V_{\Real}$ be a \textit{real} vector space of dimension 3.
Let $\SOM{2}{1}$ denote the group of unimodular linear automorphisms
of $V_{\Real}$
preserving a form $B_{\Real}$ on $V_{\Real}$ of bilinear signature $(2,1)$.
For definiteness,
we will take
\startdisp
V_{\Real}=\Real\text{-span}(\eta)\subseteq\mathrm{Lie}(\SL{2}{\Complex}),\quad
B_{\Real}=B|_{V_{\Real}},
\finishdisp
where $\beta'$ is the basis of $\mathrm{Lie}(\SL{2}{\Complex})$
defined at
\eqref{eqn:betaprime}, and $B$ is as usual the Killing form
on $\mathrm{Lie}(\SL{2}{\Complex})$.  From the fact that $\beta$
is an orthonormal set under $B$ and from \eqref{eqn:etadefn},
it is immediately verified that $B|_{\Real}$ has signature $(2,1)$.
Note also that
\startdisp
V:=V_{\Real}\otimes\Complex=\mathrm{Lie}(\SL{2}{\Complex}).
\finishdisp

By considering $\SOM{2}{1}$ as a subset of $\GL{3}{\Real}$
we obtain the \bmth\textbf{standard representation of $\SOM{2}{1}$}\ubmth.
We define $\SOM{2}{1}_{\Int}$ to be the matrices with integer
coefficients in the standard representation of $\SOM{2}{1}$.

Recall from \eqref{eqn:conjwithresptobasis} the definition of
the morphism
\startdisp
\conj_{\eta}:=\conj_{V,\eta}:
\SL{2}{\Complex}\rightarrow\SL{3}{\Real}.
\finishdisp

\begin{prop} \label{prop:gammaziso} Let $\Gamma_{\Int}$ as defined
in \eqref{defn:gammaz}.  Then the restriction of $\conj_{\eta}$
to $V_{\Real}$ provides an isomorphism
\starteqn\label{eqn:liegroupsiso}
\conj_{\eta}:\SL{2}{\Real}/\{\pm I\}\rightarrow\SOM{2}{1}^0
\finisheqn
of Lie groups.  The isomorphism of \eqref{eqn:liegroupsiso}
further restricts to an isomorphism of discrete subgroups
\starteqn\label{eqn:discretegroupsiso}
\conj_{\eta}: \conj\inv(\Gamma_{\Int})\rightarrow\SOM{2}{1}_{\Int}.
\finisheqn
As a result, $\conj_{\eta}\conj\inv$ exhibits an isomorphism
\starteqn\label{eqn:gammaintjustification}
\Gamma_{\Int}\cong\SOM{2}{1}_{\Int}.
\finisheqn
\end{prop}

The next Proposition, \ref{prop:gammaintexplicitdescrip}, is the analogue
of Proposition \ref{prop:inversematrixexplicitdescription}
for the real form of the complex group.
Proposition \ref{prop:gammaintexplicitdescrip} below
is, in contrast, almost a triviality to prove at this point,
since it can be deduced rather readily from Proposition
\ref{prop:inversematrixexplicitdescription}.

For Proposition \ref{prop:gammaintexplicitdescrip},
it is necessary to recall the group $\Xi$-subgroups of
defined in \eqref{eqn:xi12defn} and \eqref{eqn:residuematrices}.  For
each of the three $\Xi$-subgroups, we define
\starteqn\label{eqn:xi12intdefn}
(\Xi)_{\Int}=\Xi\cap\SL{2}{\Real}.
\finisheqn
The following result both justifies this notation and clarifies
the meaning of Proposition \ref{prop:gammaintexplicitdescrip}, below.
\begin{lem}
Each $(\Xi)_{\Int}$-group can be given the following description.
\starteqn\label{eqn:xiintsetdescription}
\begin{gathered}
\text{For fixed}\;
\begin{pmatrix}\overline{p}&\overline{q}\end{pmatrix},
\begin{pmatrix}\overline{r}&\overline{s}\end{pmatrix}\in
\left\{
\begin{array}{l}
\begin{pmatrix}1&1\end{pmatrix},\vspace*{0.15cm}\\
\begin{pmatrix}1&0\end{pmatrix},\vspace*{0.15cm}\\
\begin{pmatrix}0&1\end{pmatrix}
\end{array}
\right\}\subset(\SL{2}{\Int[\vecti]/(2)})^2,\\
\Xi=\red_{2}\inv\left(
\left\{\begin{pmatrix}\overline{p}&\overline{q}\\
\overline{r}&\overline{s}\end{pmatrix},\,
\begin{pmatrix}\overline{r}&\overline{s}\\
\overline{p}&\overline{q}
\end{pmatrix}\right\}\right).
\end{gathered}
\finisheqn
In order to obtain $\Xi_{12}$ in this manner,
we may take, in \eqref{eqn:xiintsetdescription},
\startdisp
\begin{pmatrix}\overline{p}&\overline{q}\end{pmatrix}=
\begin{pmatrix}1&0\end{pmatrix}\;\text{and}
\begin{pmatrix}\overline{r}&\overline{s}\end{pmatrix}
=\begin{pmatrix}0&1\end{pmatrix}
\finishdisp
Further, we may take
\startdisp
\begin{pmatrix}\overline{p}&\overline{q}\end{pmatrix}=
\begin{pmatrix}1&1\end{pmatrix},\;\text{in order to obtain $\Xi_{1}$ and $\Xi_{2}$,}
\finishdisp
and
\startdisp\begin{gathered}
\begin{pmatrix}\overline{r}&\overline{s}\end{pmatrix}
=\begin{pmatrix}0&1\end{pmatrix},\,\text{in order to obtain $\Xi_1$},\\
\begin{pmatrix}\overline{r}&\overline{s}\end{pmatrix}
=\begin{pmatrix}1&0\end{pmatrix},\,\text{in order to obtain $\Xi_2$}.
\end{gathered}
\finishdisp
\end{lem}

\begin{prop}\label{prop:gammaintexplicitdescrip}  With
$\Gamma_{\Int}$ defined as in \eqref{eqn:gammazdefn},
we have
\starteqn\label{eqn:gammaintdescription}
\conj\inv(\Gamma_{\Int})=(\Xi_{12})_{\Int}\bigcup
\frac{1}{\sqrt{2}}({\Xi_{2}})_{\Int}\begin{pmatrix}1&-1\\0&2
\end{pmatrix}.
\finisheqn
\end{prop}

From \eqref{eqn:gammaintdescription}, we deduce the analogue of
Lemma \ref{lem:indexlemma}
\begin{lem}\label{lem:indexlemmarcase}  Let $\conj\inv(\Gamma_{\Int})$
be the discrete subgroup of $\SL{2}{\Real}$ defined in \ref{eqn:gammazdefn},
and given explicitly in matrix form
in \eqref{eqn:gammaintdescription}.  All the
other notation is also as in
Proposition \ref{prop:gammaintexplicitdescrip}.
\begin{itemize}\item[(a)]  We have
\startdisp
\conj\inv(\Gamma_{\Int})\cap\SL{2}{\Int}=(\Xi_{12})_{\Int}.
\finishdisp
\item[(b)]  We have
\starteqn
[\conj\inv(\Gamma_{\Int}):(\Xi_{12})_{\Int}]=2,\quad [\SL{2}{\Int}:\Xi_{12}]=3.
\finisheqn
Explicitly, a representative of the unique non-identity right coset of
$(\Xi_{12})_{\Int}$ in $\conj\inv(\Gamma)$ is
\startdisp
\frac{1}{\sqrt{2}}
\begin{pmatrix}1&0\\1&1\end{pmatrix}\begin{pmatrix}1&-1\\0&2\end{pmatrix}=
\frac{1}{\sqrt{2}}\begin{pmatrix}1&-1\\1&1\end{pmatrix}.
\finishdisp
\end{itemize}
\end{lem}

\section{Fundamental
domain for $\SOM{2}{1}_{\Int}$ acting on $\bbH^2$
and its relation to that of $\SO{3}{\Int[\vecti]}$}
The main point of this section is that, provided the fundamental
domain $\scrG_{\Real}$ of the the standard unipotent subgroup
of $\conj\inv(\Gamma_{\Int})$ is chosen in a way that is compatible with the
choice of $\scrG$ in \eqref{eqn:scrGdefn}, then the good Grenier
fundamental domain $\scrF_{\Real}(\scrG_{\Real})$ for
$\conj\inv(\Gamma_{\Int})$ corresponding to $\scrG_{\Real}$
will have a close geometric relationship to $\scrF(\scrG)$.  Based
on the classical example of Dirichlet's fundamental domain
for $\SL{2}{\Int}$ acting on $\bbH^2$ and the Picard domain,
one might guess that we would have the equality
\starteqn
\scrF_{\Real}(\scrG_{\Real})=\scrF(\scrG)\cap\bbH^2.
\finisheqn
In fact, this intersection property cannot hold, because of the presence
of additional torsion elements (the powers of $\omega_8I_2$) in
$\conj\inv(\Gamma)$.  However, in a sense which will be made precise
in Proposition \ref{prop:funddomsrelation}, below, the next best thing holds.
Namely, the intersection
of the set consisting of \textit{two} $\Gamma$-translates of $\scrF(\scrG)$
with $\bbH^2$ equals $\scrF_{\Real}(\scrG_{\Real})$, for the choice
of $\scrG_{\Real}$ in \eqref{eqn:scrGRdefn}, below.
\vspace*{0.3cm}

In the case of $\Gamma_{\Int}\subset\mathrm{Aut}^+(\bbH^2)$, commensurable
to $\SL{2}{\Int[\vecti]}$, we have the obvious analogue
of Theorem \ref{thm:gammagoodgrenier}, defining
a good Grenier fundamental domain for the action of $\Gamma_{\Int}$.
In order to distinguish the real case
$\Gamma_{\Int}\subset\mathrm{Aut}^+(\bbH^2)$ from the complex
case, we add the subscript $\Real$ to the sets $\scrG$
$\scrF_1$, $\scrF(\scrG)$, and so write $\scrG_{\Real}$
$\scrF_{1,\Real}$, $\scrF_{\Real}(\scrG_{\Real})$.
In this case, the good Grenier fundamental domain coincides
with the classical notion of the \textit{Ford fundamental domain}
for a discrete subgroup of $\mathrm{Aut}^+(\bbH^2)$ of
finite covolume.  See, for example, \cite{iwaniec95}, p. 44.  However,
we use the terminology Grenier domain even in this context, in
order to stress the eventual connections with the higher-rank case.

\vspace*{3mm}
\noindent \bmth\textbf{Explicit Descriptions of $\scrG_{\Real}$
and $\scrF_{\Real}(\scrG_{\Real})$.}\ubmth
\begin{lem}
\begin{itemize}
\item[(a)]
We have
\startdisp
(\Gamma_{\Int})^{\phi_1}=\begin{pmatrix}1&2\Int\\
0&1\end{pmatrix}.
\finishdisp
\item[(b)]  The interval
\starteqn\label{eqn:scrGRdefn}
\scrG_{\Real}:=[0,2]
\finisheqn
is a fundamental domain for the action of $\Gamma_{\Int}^{\phi_1}$
on $\Real$ satisfying
\startdisp
\scrG_{\Real}=\overline{\intrr\scrG_{\Real}}.
\finishdisp
\item[(c)]  With $\scrG_{\Real}$ as defined in \eqref{eqn:scrGRdefn},
part (b) implies that
\starteqn\label{eqn:scrFRscrGRdescription}
\begin{aligned}
\scrF_{\Real}\left(\scrG_{\Real}\right)&=&&\{z\in\bbH^2\;
|\;0\leq x(z)\leq 2,\;y(z)^2+(x-1)^2\geq 2\}\\
&=&&\coh(\vecti,2+\vecti,\infty).
\end{aligned}
\finisheqn
\end{itemize}
\end{lem}

\vspace*{3mm}
\noindent \bmth\textbf{Geometric
relation of $\scrF_{\Real}(\scrG_{\Real})$ to $\scrF(\scrG)$.}\ubmth\;
In order to relate the fundamental domain of a subgroup
of $\SL{2}{\Real}$ acting on $\bbH^2$ to the fundamental domain of a subgroup
of $\SL{2}{\Complex}$ acting on $\bbH^3$, we consider $\bbH^2$
embedded in $\bbH^3$ as the totally geodesic surface $\bbH^2(0)$.
Note that
\startdisp
\bbH^2(0)=\{x\vecti+y\vectj\;|\; y>0\},
\finishdisp
and the actions of $\SL{2}{\Real}$ on $\bbH^2$
and $\bbH^2(0)$ are equivariant with the obvious
isomorphism
\startdisp
\bbH^2\stackrel{\cong}{\longrightarrow}\bbH^2(0),\;
\text{mapping}\;x+y\vectj\mapsto x\vecti+y\vectj.
\finishdisp
Under this isomorphism of $\SL{2}{\Real}$-homogeneous spaces,
$\scrF_{\Real}(\scrG_{\Real})$ corresponds to
\starteqn\label{eqn:scrFRscrGRHjform}
\coh(\vectj,2\vecti+\vectj,\infty)\;\text{in}\;\bbH^2(0).
\finisheqn
Because of the isomorphism, we can safely ignore the distinction
between the forms of $\scrF_{\Real}(\scrG_{\Real})$ in
\eqref{eqn:scrFRscrGRdescription} and \eqref{eqn:scrFRscrGRHjform}.

Because, as can be verified readily,
\starteqn\label{eqn:scrGRunion}
\scrG_{\Real}=\left(\scrG\cup\conj(T_1)\left(R_{\frac{\pi}{2}}^2\right)\scrG
\right)\cap
\bbH^2_{\vectj},
\finisheqn
we cannot hope that we will have the straightforward relation
\startdisp
\scrF_{\Real}\left(\scrG_{\Real}\right)=\scrF\left(\scrG\right)\cap\bbH^2_{\vectj}
\finishdisp
that we find in the classical case of  $\SL{2}{\Int[\vecti]}$ and
$\SL{2}{\Int}$.
However, we do have the next best possible relation between
the fundamental domains.

\begin{prop} \label{prop:funddomsrelation} We have the relation
\startdisp
\scrF\left(\scrG_{\Real}\right)=\left(\scrF(\scrG)\cup
\conj(T_1)\left(R_{\frac{\pi}{2}}^2\right)\scrF(\scrG)\right)\cap
\bbH^2_{\vectj}.
\finishdisp
\end{prop}

\begin{rem}  We note for possible future reference that
$\scrF_{\Real}(\scrG_{\Real})$ is the \textit{normal geodesic projection}
of the union of $\scrF(\scrG)$ and one translate
$\conj(T_1)\left(R_{\frac{\pi}{2}}^2\right)\scrF(\scrG)$ of $\scrF(\scrG)$.
This relation between the fundamental domains
is connected to the one given in Proposition \ref{prop:funddomsrelation},
though neither relation implies the other, in general.
In Figure 1, we have indicated by means of a ``right-angle" symbol
at the point $1+\sqrt{2}\vectj$ that the geodesic $\bbH^1(1+\sqrt{2}\vectj,
1+\vecti+\vectj)$ is a geodesic normal to $\bbH^2_{\vectj}$.
It would take us to far afield of our main purpose to define
the concept of \textit{normal geodesic projection} precisely,
so for the moment we restrict ourselves to mentioning that this
relation between $\scrF(\scrG)$
and $\scrF_{\Real}(\scrG_{\Real})$ may be of some
use in relating spectral expansions in the complex case
to spectral expansions in the real case.
\end{rem}

\section{Spectral Zeta Functions}
This section discusses a potential application of the results of
the paper and indicates a future line of investigation building
on this work.  Jorgenson and Lang, in works such as \cite{jol01},
\cite{jol05} (see the introduction to the latter work especially),
and \cite{jol06}, have laid out and begun
to pursue an ambitious program of using heat kernel
analysis to associate additive spectral zeta functions
to quotients of symmetric spaces.  When completed, this
theory would subsume the basic theory of the Riemann zeta
function and Selberg zeta function (among others), and clarify
the relationship between the zeta functions arising at different
geometric levels.  The main component of the program
is obtaining a theta inversion formula.

In \cite{jol05}, which carries out the derivation of the theta
inversion formula for the special case of
\startdisp
X=\Gamma\backslash
G/K=\SL{2}{\Int[\vecti]}\backslash\SL{2}{\Complex}/\mathrm{SU}(2,\Complex),
\finishdisp
the authors compute the regularized trace of an integral
operator on functions on $X$.  The kernel of the integral
operator is $\mathbf{K}_{t,X}(z,w)$, the heat kernel on $X$.
The trace of such an integral operator is defined to be the integral
on the diagonal
\startdisp
\int_{X}\mathbf{K}_{t,X}(z,z)\intd z.
\finishdisp
Although this integral is infinite, because of the cusp of $X$,
the integrals over sets $X_Y$ approximating by covering $X$ only
up to some fixed finite ``distance" in the cusp are finite
and diverge logarithmically in $Y$.  That is,
\starteqn\label{eqn:truncationlimit}
\lim_{Y\rightarrow\infty}
\int_{X_Y}\mathbf{K}_{t,X}(z,z)\intd z-c_1(t)\log Y\;\text{exists
as a $\Complex$-valued function of $t$.}
\finisheqn
where $c_1(t)$ is a factor, constant in $Y$, and determined in \cite{jol05}.
For the purposes of such an integration, we can
replace $X$ with a suitable fundamental domain $\scrF$.
The fundamental domain $\scrF$ is
an analytic model of $X$ in its universal covering
space $\bbH^3$--see \S\ref{sec:explicitfd}, below, for a
precise definition of fundamental domain.  Similarly,
we replace the truncation $X_Y$ of $X$ with a
matching truncation $\scrF_Y$ of $\scrF$.

To obtain the theta inversion formula, the
limit of \eqref{eqn:truncationlimit} is computed
in two ways.  One computation is
from the expression of the heat kernel as the periodized
heat kernel on the universal covering space $\bbH^3$.
This method of computing the limit in \eqref{eqn:truncationlimit} yields
\starteqn\label{eqn:invertedtheta}
e^{-2t}(4t)^{-\half}\Theta^{\rm NC}(1/t)+\Theta^{\rm Cus}(1/t).
\finisheqn
In \eqref{eqn:invertedtheta}, $\Theta^{\rm NC}(1/t)$, the
inverted theta series, is
defined in terms of invariants
of certain $\Gamma$-conjugacy classes in $\Gamma$, while
$\Theta^{\rm Cus}(1/t)$, the inverted theta
integral, is a sum of products composed
of special values of $\zeta_{\Rational(\vecti)}$, constants
similar to Euler's $\gamma$, and single integrals whose
Gauss transforms are exact.  (We refer to
\newcounter{ranrom3}\setcounter{ranrom3}{14}\S\Roman{ranrom3}.7,
of \cite{jol05}, for exact definitions of
$\Theta^{\rm NC}(1/t)$, $\Theta^{\rm Cus}(1/t)$ and the other
terms in the theta relation.)  The other method of computing
the limit in \eqref{eqn:truncationlimit}
is from the expansion of $\mathbf{K}_{t,X}(z,z)\intd z$
in terms of the spectrum of the Laplacian $\BDelta_{X}$.
This second method of computing the limit
of \eqref{eqn:truncationlimit} yields
\starteqn\label{eqn:theta}
\theta_{\rm Cus}(t)+1+\theta_{\rm Eis}(t),
\finisheqn
where $\theta_{\rm Cus}(t)$ is the theta series
$\sum_{j=1}^{\infty}e^{-\lambda_jt}$ and $\lambda_j$
are the eigenvalues of $\BDelta_{X}$, and $\theta_{\rm Eis}(t)$
is what remains as the limit of the integral of the convolution
of $\mathbf{K}_{t,X}(z,w)$ with certain Eisenstein series, once
the term $c_1(t)\log(Y)$ has been subtracted.
Setting equal the two expressions, \eqref{eqn:invertedtheta}
and \eqref{eqn:theta}, for the same limit \eqref{eqn:truncationlimit},
we obtain the theta inversion formula for $X$,
\starteqn\label{eqn:thetainversionforX}
e^{-2t}(4t)^{-\half}\Theta^{\rm NC}(1/t)+\Theta^{\rm Cus}(1/t)=
\theta_{\rm Cus}(t)+1+\theta_{\rm Eis}(t).
\finisheqn
Next, note that there is an infinite sequence of arithmetic
quotients
\startdisp
X_n=\SL{n}{\Int[\vecti]}\backslash\SL{n}{\Complex}/\mathrm{SU}(n),\;
n>1,
\finishdisp
having $X=X_2$ as its first nontrivial member.  Generalizations
of \eqref{eqn:thetainversionforX} to $X_n$ for $n>2$
are discussed in \cite{jol06}.  In order to obtain
exact formulas analogous to \eqref{eqn:thetainversionforX},
we would have to integrate over a fundamental domain, rather
than over an appoximating Siegel set, which is a more common
analytic model in the literature.

In the present work, we initiate an extension of the Jorgenson-Lang
project to the sequence of arithmetic quotients
\starteqn\label{eqn:ourlocallysymmetricspace}
\SO{n}{\Int[\vecti]}\backslash \SO{n}{\Complex}/\SOR{n}
\finisheqn
and related arithmetic quotients of real forms of the symmetric
spaces.  The main results of the present paper are restricted to the
group theory (Propositions \ref{prop:inversematrixexplicitdescription} and
\ref{prop:gammaintexplicitdescrip})
and fundamental domains
(Propositions \ref{prop:gammafunddom} and \ref{prop:funddomsrelation})
in the first case
of $n=2$.  Nevertheless, some of the intermediate results
are couched in a more general terminology and notation,
with a view towards building upwards from the case $n=2$, to the
case of a general $n$.  Thus, our project
includes a natural extension and generalization
of Grenier's work in \cite{grenier88}
and \cite{grenier93} to a different sequence of symmetric spaces.

 The identification
\startdisp
\SL{2}{\Complex}/\{\pm I\}
\stackrel{\cong}{\longrightarrow} \SO{3}{\Complex},
\finishdisp
allows us to view the theta inversion relation (conjecturally)
associated with the case $n=2$ in \eqref{eqn:ourlocallysymmetricspace}
as a theta inversion relation associated with a quotient
of $\SL{2}{\Complex}/K$ by an arithmetic subgroup different
from, but still commensurable, the ``standard" arithmetic
subgroup $\SL{2}{\Int[\vecti]}$.  The results
of this paper will, it is hoped, enable future investigations to apply
the machinery developed in \cite{jol05} to this ``nonstandard" arithmetic
subgroup $\conj\inv(\SO{3}{\Int[\vecti]})$ of
$\SL{2}{\Complex}$ to obtain the corresponding
theta function.

\begin{wrapfigure}{r}{5cm}
\includegraphics[scale=0.5]{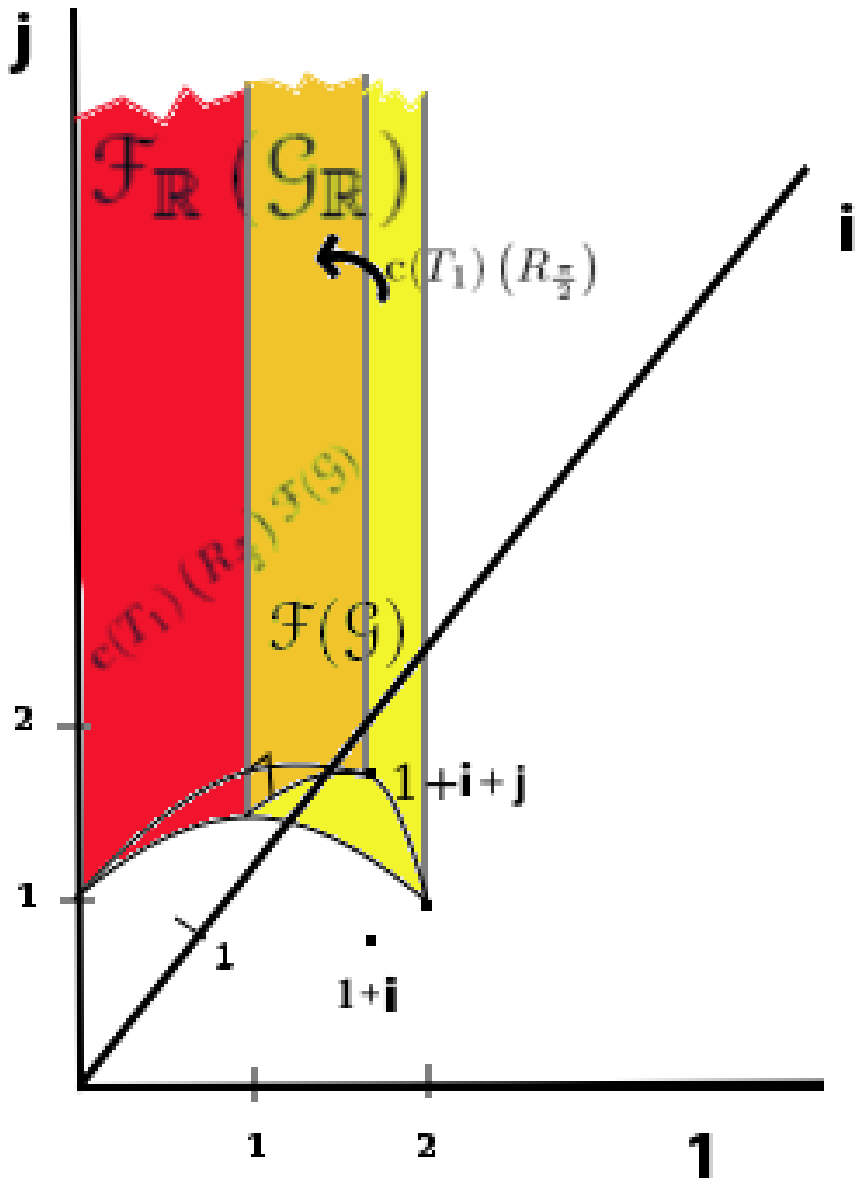}
   \caption{Fundamental domains for $\Gamma$ and $\Gamma_{\Int}$,
   with illustration of how $\conj(T_1)\left(R_{\frac{\pi}{2}}^2\right))$
   rotates the subset $\scrF(\scrG)\cap\bbH^2$ of
   $\scrF_{\Real}(\scrG_{\Real})$ into the other half of
   $\scrF_{\Real}(\scrG_{\Real})$.}
   \end{wrapfigure}

\begin{singlespace}
\bibliographystyle{amsalpha}
\bibliography{../memoirsbib}
\addcontentsline{toc}{chapter}{Bibliography}
\end{singlespace}
\end{document}